\documentclass[12pt,a4paper,reqno]{amsart}
\usepackage{amsmath}
\usepackage{amsfonts}
\usepackage{amssymb}
\usepackage{bm} 

\newcommand\R{{\mathbf{R}}}
\newcommand\C{{\mathbf{C}}}
\newcommand\Z{{\mathbf{Z}}}
\newcommand\K{{\mathcal{K}}}
\newcommand\eps{{\varepsilon}}
\newcommand\hi{{\operatorname{hi}}}
\newcommand\rad{{\operatorname{rad}}}
\newcommand\med{{\operatorname{med}}}
\newcommand\dist{{\operatorname{dist}}}
\newcommand\lo{{\operatorname{lo}}}
\newcommand\wlim{{\mathop{\operatorname{w-lim}}_{t \to +\infty}}\ }
\newcommand\Wlim{{\mathop{\operatorname{w-lim}}_{T \to +\infty}}\ }

\theoremstyle{plain}
  \newtheorem{theorem}[subsection]{Theorem}
  
  \newtheorem{proposition}[subsection]{Proposition}
  \newtheorem{lemma}[subsection]{Lemma}
  \newtheorem{corollary}[subsection]{Corollary}

\theoremstyle{remark}
  \newtheorem{remark}[subsection]{Remark}
  
  \newtheorem{example}[subsection]{Example}

\theoremstyle{definition}
  \newtheorem{definition}[subsection]{Definition}

\include{psfig}

\begin{document}

\title[A compact attractor for high-dimensional NLS]{A (concentration-)compact attractor for high-dimensional non-linear Schr\"odinger equations}
\author{Terence Tao}
\address{Department of Mathematics, UCLA, Los Angeles CA 90095-1555}
\email{tao@@math.ucla.edu}
\subjclass{35Q55}

\vspace{-0.3in}
\begin{abstract}
We study the asymptotic behavior of large data solutions to Schr\"odinger equations
$i u_t + \Delta u = F(u)$ in $\R^d$,
assuming globally bounded $H^1_x(\R^d)$ norm (i.e. no blowup in the energy space), in high dimensions $d \geq 5$ and with nonlinearity which is energy-subcritical and mass-supercritical.  In the spherically symmetric case, we show that as $t \to +\infty$, these solutions split into a radiation term that evolves according to the linear Schr\"odinger equation, and a remainder which converges in $H^1_x(\R^d)$ to a compact attractor, which consists of the union of spherically symmetric almost periodic orbits of the NLS flow in $H^1_x(\R^d)$.  This is despite the total lack of any dissipation in the equation.
This statement can be viewed as weak form of the ``soliton resolution conjecture''.  We also obtain a more complicated analogue of this result for the non-spherically-symmetric case.  As a corollary we obtain the ``petite conjecture'' of Soffer in the high dimensional non-critical case.
\end{abstract}

\maketitle

\section{Introduction}

The purpose of this paper is to establish some asymptotic properties of bounded-energy solutions of non-linear Schr\"odinger (NLS) equations 
\begin{equation}\label{nls}
i u_t + \Delta u = F(u)
\end{equation}
with moderate (but possibly focusing) nonlinearity $F: \C \to \C$ and high dimension $d$; we allow the nonlinearity $F$ to be focusing in nature as long as the energy remains bounded.  The main result is to describe a certain \emph{compact attractor} for the NLS flow, although the definition of ``attractor'' (and ``compact'') needs to be defined properly in this dispersive (and translation-invariant) context.

\subsection{Assumptions on the equation}

We shall only consider NLS equations \eqref{nls} which obey the following hypotheses:

\begin{itemize}

\item (High dimension) $d \geq 5$.

\item (Hamiltonian structure) There exists a $C^1$ function $G: \R^+ \to \R$ with $G(0)=0$ such that $F(z) = G'(|z|^2) z$ for all $z \in \C$.

\item (Power-type nonlinearity) There exists an exponent $p > 1$, a constant $C_0>0$ and a H\"older regularity index $0 < \theta \leq \min(p-1,1)$ for which we have the estimates
\begin{align}
|F(z)| &\leq C_0 |z|^p  \label{fpower}\\
|F'(z)| &\leq C_0 |z|^{p-1}  \label{fpower-2}\\
|F'(z) - F'(w)| &\leq C_0  |z-w|^\theta (|z| + |w|)^{p-1-\theta} \label{fpower-3}
\end{align}
for all $z,w \in \C$.  Here we view the differential $F'(z)$ of $F$ at $w$ as a real-linear map from $\C$ to $\C$.  

\item (Mass-supercriticality) We have $p > 1 + \frac{4}{d}$.

\item (Energy-subcriticality) We have $p < 1 + \frac{4}{d-2}$.

\end{itemize}

{\bf Important convention.} Throughout this paper we fix $d,p,\theta,C_0,F$ and we shall always assume the above hypotheses to be in effect.  Also, all quantities in this paper are implicitly assumed to depend on the dimension $d$, the exponent $p$, the H\"older regularity $\theta$, and the constant $C_0$.

Note that we have made no assumptions about the sign of the nonlinearity $F$ or the potential function $G$.  Important examples of NLS of the above type to keep in mind are:

\begin{itemize}

\item (Quadratic case) $d=5$ and $p=2$.

\item (Coercive case) $\liminf_{x \to +\infty} G(x) / x^{(d+2)/d} \geq 0$.

\item (Defocusing case) $F(z) = +|z|^{p-1} z$. (This is coercive.)

\item (Focusing case) $F(z) = -|z|^{p-1} z$.  (This is non-coercive.)

\end{itemize}

We also note the \emph{completely integrable case} when $d=1, p=3$, and the nonlinearity is either focusing or defocusing; this case is not, strictly speaking, covered by the above hypotheses (the dimension is too low), but is better understood than most other NLS equations and serves as motivation for the soliton resolution conjecture which we discuss later.

The NLS equation \eqref{nls} is manifestly translation invariant.  The assumption of Hamiltonian structure also gives us the symmetries of \emph{phase invariance} $u \mapsto e^{i\alpha} u$ and \emph{Galilean invariance}
\begin{equation}\label{gal}
 u(t,x) \mapsto e^{i v \cdot x / 2} e^{-i |v|^2 t / 4} u(t, x-vt)
 \end{equation}
for any $v \in \R$.  The Hamiltonian structure also gives several conserved quantities, including the
\emph{mass}
\begin{equation}\label{mass-def}
M(u) := \int_{\R^d} |u(t,x)|^2\ dx
\end{equation}
and the \emph{Hamiltonian}
\begin{equation}\label{hamil}
 H(u) := \int_{\R^d} \frac{1}{2} |\nabla u(t,x)|^2 + \frac{1}{2} G( |u(t,x)|^2 )\ dx.
\end{equation}

Throughout this paper we shall be working in the \emph{energy space} $H := H^1_x(\R^d \to \C)$, which is a Hilbert space with inner product
$$ \langle u, v \rangle_H := \int_{\R^d} u(x) \overline{v(x)} + \nabla u(x) \cdot \overline{\nabla v(x)}\ dx.$$

\begin{remark}\label{energy-1}
From Sobolev embedding and the energy-subcritical nature of $p$ we observe that the mass $M(u)$ and Hamiltonian $H(u)$ are finite for any $u \in H$.  Conversely, if we assume the NLS equation is coercive, then a standard application of the Gagliardo-Nirenberg inequality reveals that any function $u$ with finite mass and Hamiltonian lies in $H$.  However, in non-coercive cases, such as the focusing case, it is certainly possible for the $H$ norm to blow up in finite time even with finite mass and Hamiltonian, see \cite{glassey}.
\end{remark}

\subsection{Assumptions on the solution}

We shall only consider \emph{bounded-energy} solutions to \eqref{nls}, although we allow this energy bound to be arbitrarily large.  More precisely, we have

\begin{definition}[Solutions]  A \emph{bounded-energy strong solution to \eqref{nls}}, or \emph{solution} for short, will be any function $u \in C^0_t H^1_x(I \times \R^d)$ on a non-empty time interval $I \subset \R$ taking values continuously in the energy space $H$ such that
$$ u(t_1) = e^{i(t_1-t_0)\Delta} u(t_0) - i \int_{t_0}^{t_1} F(u(t))\ dt$$
for all $t_0, t_1 \in I$, where we of course adopt the convention that $\int_{t_0}^{t_1} = -\int_{t_1}^{t_0}$ if $t_1 < t_0$, and such that the \emph{energy}
\begin{equation}\label{energy-bound}
E(u) := \sup_{t \in I} \| u(t) \|_H^2
\end{equation}
is finite.
Here $e^{it\Delta}$ is the free Schr\"odinger propagator, defined via the Fourier transform
$$
\hat f(\xi) := \int_{\R^d} e^{-ix \cdot \xi} f(x)\ dx
$$
by
\begin{equation}\label{prop-fourier}
\widehat{e^{it\Delta} f}(\xi) := e^{-it|\xi|^2} \hat f(\xi)
\end{equation}
or more directly as
\begin{equation}\label{prop}
e^{it\Delta} f(x) = \frac{1}{(4\pi i t)^{d/2}} \int_{\R^d} e^{i|x-y|^2/4t} f(y)\ dy.
\end{equation}
We say that a solution is \emph{forward-global} if $I$ contains $[0,+\infty)$, and \emph{global} if  $I = \R$.
\end{definition}

\begin{remark} One can of course talk about backward-global solutions, which contain $(-\infty,0]$, but because of the
time-reversal symmetry $u(t,x) \mapsto \overline{u(-t,x)}$ all the results here for forward-global solutions immediately have counterparts for backward-global solutions (and hence, by concatenation, for global solutions).
\end{remark}

We make the trivial observation that the restriction of any solution to a sub-interval is still a solution.  Also observe that the propagators $e^{it\Delta}$ are unitary on $H$.

\begin{remark}  Note that a solution $u: \R \to H$ whose $H$ norm goes to infinity as $t \to \pm \infty$ would not be considered a global solution in our notation, though it is a solution on any compact sub-interval of $\R$.
\end{remark}

We list some basic facts about solutions below.  Define an exponent pair $(q,r)$ to be \emph{admissible} if $\frac{2}{q} + \frac{d}{r} = \frac{d}{2}$ and $2 \leq q,r \leq \infty$.  

\begin{theorem}[Local existence and uniqueness]\label{local} \ \ 
\begin{itemize}
\item (Local existence) If $t_0 \in \R$ and $B \subset H$ is bounded, then there exists an open time interval $I$ containing $t_0$ such that for every $u_0 \in B$ there exists a solution $u: I \to H$ such that $u(t_0)=u_0$.  Furthermore the map $u_0 \mapsto u(t)$ is Lipschitz continuous on $B$ for all $t \in I$.  If $B$ is a sufficiently small neighbourhood of the origin, one can take $I=\R$.
\item (Uniqueness) If two solutions $u: I \to H$, $\tilde u: I \to H$ agree on at least one time, then they are equal for all time.
\item (Strichartz regularity) If $u: I \to H$ is a solution, $J$ is a compact sub-interval of $I$, and $(q,r)$ is an admissible pair, then $u, \nabla u \in L^q_t L^r_x(J \times \R^d)$.
\item (Finite time blowup condition) If a solution $u: I \to H$ with finite future endpoint $T_+ := \sup I < +\infty$ cannot be extended beyond $T_+$, then $\|u(t)\|_H$ goes to infinity as $t \to T_+$ from below.  Similarly for solutions which cannot be extended beyond their finite past endpoint $T_- = \inf I > -\infty$.
\item (Conservation laws) The mass $M(u(t))$ and Hamiltonian $H(u(t))$ are constant for $t \in I$.
\end{itemize}
\end{theorem}

\begin{proof} See \cite{gv1}, \cite{gv:localreference}, \cite{cwI}, \cite{caz}, or \cite{taobook} for the local existence, regularity, conservation laws, and finite time blowup (to obtain the endpoint regularity $(q,r) = (2, \frac{2d}{d-2})$, one needs the endpoint Strichartz estimate in \cite{tao:keel}).  The uniqueness claim is proven in \cite{katounique}.
\end{proof}

\begin{remark} The above theorem, combined with Remark \ref{energy-1}, shows that in the coercive case that for any $t_0 \in \R$ and $u_0 \in H$ there is a unique global solution $u: \R \to H$ with $u(t_0) = u_0$; see e.g. \cite{gv:localreference}.  But in non-coercive cases with large initial data, global solutions need not exist; see \cite{glassey}, \cite{ogawa}.  
\end{remark}

We define the non-linear flow maps $S(t)$ on $H$ for $t \in \R$ by setting $S(t) u(0) := u(t)$ whenever
$u: [0,t] \times \R^d \to \C$ is a solution.  These maps are not necessarily globally defined (except in the coercive case), but from the above theorem we see that they are continuous and obey the group law $S(t) S(t') = S(t+t')$ on their domain of definition, and for any bounded set in $H$ the $S(t)$ are defined for all sufficiently small $t$.

\subsection{The soliton resolution conjecture}

Suppose we have a forward-global solution $u: I \to H$ to the NLS \eqref{nls}.  A natural question then
arises as to what the asymptotic behaviour of $u(t)$ is as $t \to +\infty$.  In the case of small energy or defocusing nonlinearity, the answer is known:

\begin{theorem}[Scattering]\label{scat}  Let $u_0 \in H$ and $t_0 \in \R$, and assume either that $\|u_0\|_H$ is sufficiently small, or that the nonlinearity is defocusing.  Then there is a unique global solution $u: \R \to H$ with $u(t_0) = u_0$.  Furthermore there is a unique \emph{scattering state} $u_+ \in H$ such that
$$ \lim_{t \to +\infty} \| u(t) - e^{it\Delta} u_+ \|_H = 0.$$
\end{theorem}

\begin{proof} For the small data case, see \cite{gv:localreference}, \cite{strauss}, \cite{caz}, or \cite{taobook}.  The defocusing case was proven in \cite{gv:scatter} (see also \cite{borg:scatter}, \cite{ckstt:scatter}, \cite{visan:scatter}, \cite{taobook}).  The proofs of the defocusing result rely on Morawetz inequalities, which do not have a favourable sign in other cases, including the focusing case and even some coercive cases.
\end{proof}

Equivalently, we may write
\begin{equation}\label{scatter} u(t) = S(t-t_0) u_0 = e^{it\Delta} u_+ + o_H(1)
\end{equation}
where we use $o_H(1)$ to denote a time-dependent function which goes to zero in $H$ norm as $t \to +\infty$.

In the focusing case and with large data $u_0$, the above theorem fails for at least two reasons.  Firstly, as mentioned earlier, the solution can blow up in finite time, especially if the Hamiltonian
is negative; see \cite{glassey}, \cite{ogawa}.  Secondly, even if the solution remains global (or at least forward-global), it does not necessary scatter to a free solution $e^{it\Delta} u_+$.  This can be seen by considering \emph{stationary soliton solutions} of the form $u(t,x) = Q(x) e^{i\omega t}$, where $\omega > 0$ is a constant and $Q \in H$ solves the elliptic  equation
\begin{equation}\label{gse}
 \Delta Q + |Q|^{p-1} Q = \omega Q.
\end{equation}
(For a construction of such solutions, see \cite{blions}, \cite{wein1}.)  One can also apply Galilean symmetry \eqref{gal} to create travelling soliton solutions, which at time $t$ would be localised near $x_0 + vt$ for some $x_0, v \in \R^d$.  Furthermore, it is possible to create \emph{multisoliton} solutions which as $t \to +\infty$ resemble superpositions of $J$ divergent traveling solitons for any given $J \geq 1$; see \cite{perelman}, \cite{schlag}, \cite{tsai-2} for some constructions of such solutions for various choices of $d,p,F$.  Finally, it is possible in some cases to superimpose a free solution $e^{it\Delta} u_+$ with a soliton or multisoliton solution, at least for sufficiently late times.

It is tentatively conjectured that the above behaviour is in fact generic.  This leads to an (imprecise) \emph{soliton resolution conjecture}, that for ``generic'' large global solutions, 
the evolution asymptotically decouples into the superposition of divergent solitons, a free radiation
term $e^{it\Delta} u_+$, and an error which goes to zero at infinity (cf. \eqref{scatter}).  We leave questions such as the regularity and decay class of the solution, the sense in which the error goes to zero, and the definition of ``generic'' as deliberately vague.  Indeed, our understanding of this conjecture is still very poor (even with strong additional assumptions such as spherical symmetry and coercive nonlinearity), with the majority of results being concentrated either on the small data or defocusing cases (in which no solitons are present), or when the solution is very close to a soliton or multisoliton solution, especially if the solitons are generated by a ground state.  See \cite{soffer-icm} for some further discussion of this conjecture (referred to there as the \emph{grand conjecture}); see also \cite{tao:compact}.

As just mentioned, there is little direct progress on the soliton resolution conjecture for generic large data (not close to any soliton or multisoliton).  However one can consider weakening the conjecture by asking instead for an asymptotic resolution into a free solution $e^{it\Delta} u_+$, an error, and some sort of ``pseudo-multisoliton'' which exhibits behaviour similar to that of a multisoliton.  This type of conjecture is easiest to formalise in the case of spherically symmetric solutions, in which travelling solitons are precluded and the only multisoliton which is expected to be relevant is that of a single soliton placed at the origin.  But in principle one could also imagine multiple solitons of different amplitudes and widths all superimposed on each other at the origin, or more generally some sort of exotic ``breather'' solution which is periodic or almost periodic, but which does not have the explicit form $Q(x) e^{i\omega t}$.  While such solutions are expected to be very unstable, and in fact probably do not exist for most nonlinearities, we do not know how to rule them out with present technology.  Thus we can try to weaken the conjecture in this case by allowing the pseudo-soliton component to merely be almost periodic in time, rather than be an actual soliton.  As we shall see, this weakened statement is related to the \emph{petite conjecture} in \cite{soffer-icm}.

\subsection{Main results in the radial case}  

Our first set of results (which we prove in Section \ref{radial-sec}) establishes the petite conjecture in the spherically symmetric case, by showing the existence of a compact attractor for the non-radiating component of the evolution.  More precisely, we have

\begin{theorem}[Compact attractor, spherically symmetric case]\label{attract-rad} Let $E > 0$.  Then there exists a compact subset $\K_{E,\rad} \subset H$ which is invariant under the NLS flow (thus $S(t)$ is well-defined and is a homeomorphism on $\K_{E,\rad}$ for all $t \in \R$), and such that for every
spherically symmetric forward-global solution $u$ of energy at most $E$, there exists a unique \emph{radiation state} $u_+ \in H$ such that
\begin{equation}\label{tsim}
\lim_{t \to +\infty} \dist_{H} ( u(t) - e^{it\Delta} u_+, \K_{E,\rad} ) = 0.
\end{equation}
Here and in the sequel we write $\dist_H( f, K ) := \inf \{ \|f-g\|_H: g \in K \}$ for the distance between $f$ and $K$.
\end{theorem}

Thus $\K_{E,\rad}$ is a compact attractor for spherically symmetric solutions of energy at most $E$, once the effect of the radiation term $e^{it\Delta} u_+$ is removed\footnote{It is essential that we remove radiation, otherwise the concept of a compact attractor is incompatible with the time-reversibility of the NLS equation.}.
In other words, for spherically symmetric forward-global solutions $u$ of energy at most $E$, we have a decomposition\footnote{This decomposition can be regarded as a nonlinear analogue of the spectral decomposition of a linear Schr\"odinger operator with potential into continuous (dispersive) and pure point (almost periodic) components (cf. the RAGE theorem).  With this perspective, the point of the high dimension hypothesis $d \geq 5$ is to rule out ``nonlinear resonances''.} of the form
\begin{equation}\label{uw-decomp}
 u(t) = e^{it\Delta} u_+ + w(t) + o_{H}(1)
\end{equation}
where $w(t)$ ranges in the fixed compact set $\K_{E,\rad}$ for all times $t$.  Note that we do not assert that $w$ itself evolves by NLS (which would make $w$ an almost periodic solution); the problem is that the radiation terms may cause significant long-term drift in the ``secular modes'' or ``modulation parameters'' of the compact attractor $\K_{E,\rad}$.  In high dimension one expects that the strong dispersive properties of the equation will in fact rule out this scenario, but we were unable to do so here (it seems to require a linearised stability analysis of the almost periodic solutions, which we do not know how to do).

\begin{remark} A significantly weaker variant of the above theorems in \cite{tao:compact}, in the case of focusing NLS with $p=3, d=3$ (which is not covered in the analysis here).  In that paper, an attractor $\K$ was constructed in the $\dot H^1_x(\R^3)$ topology rather than the $H^1_x(\R^3)$ topology.  Also, the solutions in $\K$ were known to be uniformly smooth and enjoy some weak uniform decay at infinity, but were not known to be almost periodic, and $\K$ was not known to be compact in $H^1_x(\R^3)$.
\end{remark}

\begin{remark}\label{splash} Note that the theorem provides no information about the \emph{rate} of convergence to the compact attractor.  Indeed we expect this rate of convergence to be highly non-uniform, depending in a discontinous way on the initial data $u(0)$.  To give an example in the focusing case, suppose $u(0)$ was equal to the ground state $Q$ (which is known to be orbitally unstable for the range of exponents $p$ under consideration, see e.g. \cite{shatah}).  Then $u(t)$ will lie in the circle $\{ e^{i\alpha} Q: \alpha \in \R \}$, which we have already observed to lie in $\K_{E,\rad}$.  If however we perturb the initial data $u(0)$ by a small amount (of size $\eps$ in the $H$ norm, say), then a typical scenario would then be that after a relatively long time (e.g. of time $\log \frac{1}{\eps}$, or perhaps $\eps^{-c}$ for some $c > 0$) the solution would eventually move away from this circle, and would most likely collapse entirely into radiation.  Thus we see that the time required to reach the asymptotic state can be arbitrarily large as $\eps \to 0$, leading to a discontinuity in the decay rates in Theorem \ref{attract-rad}.  In particular, we do \emph{not} expect the compact attractor $\K_{E,\rad}$ to be orbitally stable.
\end{remark}

As one consequence of the above theorems we obtain the \emph{petite conjecture} of Soffer \cite{soffer-icm} in the radial high-dimensional case:

\begin{definition}[Almost periodic solutions]  A solution $u: I \times \R^d \to \C$ is \emph{almost periodic} if the orbit $\{ u(t): t \in I \}$ is precompact in $H$.  (See Appendix \ref{appendix} for further discussion of precompact sets in $H$.)
\end{definition}

\begin{example} The global soliton solution $u(t,x) = Q(x) e^{it}$ to a focusing NLS is almost periodic, as is any translate or rescaling of this soliton solution.  Any other hypothetical periodic or quasiperiodic ``breather'' solution to an NLS would also qualify as being almost periodic.
If one applies a Galilean transformation to give such solitons or breathers a non-zero velocity, then the solution is no longer almost periodic.  Since the set $\K_{E,\rad}$ is compact and invariant, any initial data $u(0)$ in $\K_{E,\rad}$ gives rise to a global almost periodic solution.  Thus if we could demonstrate that the only spherically symmetric almost periodic solutions were soliton solutions, Theorem \ref{attract-rad} would yield the soliton resolution conjecture in the spherically symmetric case.
\end{example}

\begin{corollary}[Petite conjecture, radial case]\label{petite-rad} Let $u$ be a spherically symmetric forward-global solution, and let $u_+$ be the radiation state.  Then the following are equivalent:
\begin{itemize}
\item[(i)] $u$ is almost periodic.
\item[(ii)] $u_+ = 0$ (i.e. $u$ is future non-radiating).
\item[(iii)] $u$ is spatially localised in the sense that\footnote{We urge the reader to pay careful attention to the order of limits in this paper, as these orderings will play a crucial role in our results and arguments.  For instance, the statement here has trivial content if the two limits are reversed.}
$$ \lim_{R \to +\infty} \limsup_{t \to +\infty} \int_{|x| > R} |u(t,x)|^2\ dx = 0.$$
\item[(iv)] $u$ is spatially localised in the sense that
$$ \lim_{R \to +\infty} \limsup_{t \to +\infty} \int_{|x| > R} |u(t,x)|^2 + |\nabla u(t,x)|^2\ dx = 0.$$
\end{itemize}
\end{corollary}

As remarked in \cite{soffer-icm}, this type of result can be regarded as a nonlinear analogue of the RAGE theorem relating dispersion and bound states for linear Schr\"odinger equations with potential.  We shall prove this Corollary in Section \ref{radial-sec}.

\begin{remark} We have not specified exactly what the compact attractor $\K_{E,\rad}$ is (although our arguments in principle provide an explicit description).  Indeed it might not be unique.  However, a simple compactness argument shows that there is a unique \emph{minimal} choice\footnote{Indeed, this minimal choice is simply the closure of the collection of all limit points $\lim_{n \to \infty} u(t_n) - e^{it_n\Delta} u_+$ for forward-global solutions $u$ of energy at most $E$, though this is clearly an unsatisfactory characterisation for $\K_{E,\rad}$.} for $\K_{E,\rad}$.  As remarked earlier, every element of $\K_{E,\rad}$ gives rise to a global almost periodic spherically symmetric solution of energy at most $E$ .  In the converse direction, by Corollary \ref{petite-rad} and Theorem \ref{attract-rad}, any limit point $\lim_{n \to \infty} u(t_n)$ of a forward-global almost periodic spherically symmetric solution $u: \R^+ \times \R^d \to \C$ (where $t_n$ is a sequence of times going to infinity) must lie in $\K_{E,\rad}$.  But these two observations do not fully pin down what $\K_{E,\rad}$ is; for instance, if the unstable manifold for the orbit of one soliton intersects the stable manifold for another, it is not clear to the author whether the intersection of these two manifolds necessarily lies in $\K_{E,\rad}$ or not, the problem being that there may be solutions which exhibit arbitrary amounts of ``Arnold diffusion'' back and forth between the two soliton orbits.  
\end{remark}

\begin{remark}
From Theorem \ref{scat} we can at least say that the zero solution $0$ is an isolated point in the minimal $\K_{E,\rad}$.  In light of the results in \cite{merlekenig}, it is also likely that after the zero solution, the next closest element of $\K_{E,\rad}$ to the origin (if it exists at all) arises from a ground state $Q$.  The instability of the ground state should then imply that the circle $\{ e^{i\alpha} Q: \alpha \in \R \}$ is an isolated connected component of $\K_{E,\rad}$, after restricting to the surface cut out by the mass and Hamiltonian conservation laws, and possibly after assuming some sign conditions on the nonlinearity.  We will however not formalise these assertions here.
\end{remark}

\begin{remark} In Theorem \ref{scat}, there is a homeomorphism between the initial data $u_0$ and the radiation state $u_+$.  We do not expect this type of correspondence to persist in the large data case, in the presence of bound states.  Firstly, the discussion in Remark \ref{splash} suggests that the map $u_0 \to u_+$ is likely to contain discontinuities, for instance at the ground state $Q$.  Secondly, there is the (somewhat strange) possibility that two initial data $u_0, u'_0$ might lead to forward-global solutions $u, u'$ which are asymptotically equivalent in the sense that $u(t)-u'(t)$ converges to zero in $H$ norm as $t \to +\infty$.  While such a scenario seems unlikely, the author was unable to argue (even heuristically) why it could not occur (although it does not seem to be possible in the completely integrable case $d=1,p=3$).  
\end{remark}

\begin{remark} Our methods actually give an explicit rate of decay for the spatial localisation, thus if $u$ is an almost-periodic forward-global spherically symmetric solution of energy at most $E$ then we have
$$ \limsup_{t \to +\infty} \int_{|x| > R} |u(t,x)|^2 + |\nabla u(t,x)|^2\ dx \leq c_E(R)$$
for some explicit quantity $c_E(R)$ which goes to zero as $R \to \infty$.  It would be of interest to obtain a good bound for this rate of decay, such as a polynomial decay $R^{-\eps}$.  Based on the observation that solitons are rapidly decreasing in space, one might even hope to get much more rapid decay, i.e. $O_N(R^{-N})$ for all $N > 0$.
Unfortunately our methods here only give a much weaker decay, something like $1/\log^c R$.  One important milestone might be to obtain a decay better than $1/R^2$ for the mass density $|u(t,x)|^2$, as this would then place the weakly bound component of the solution in the scattering space $\Sigma = \{ u: xu, u, \nabla u \in L^2_x(\R^d) \}$ and allow for tools such as the pseudoconformal identity to be applied.
\end{remark}

\subsection{Main results in the general case}  We now turn to the general case, in which no spherical symmetry is assumed.  The key difficulty here is that the class of solutions of energy $E$ is now translation-invariant, and so the notion of almost periodicity needs to be replaced by a more general notion which is both translation-invariant and also closed under certain ``superposition'' operations.

\begin{definition}[Symmetry group] Given any $h \in \R^d$, we let $\tau_h:H \to H$ be the (unitary) shift operator $\tau_h f(x) := f(x-h)$, and we let $G := \{ \tau_h: h \in \R^d \}$ be the associated translation group.  Note that this is a non-compact Lie group and so there is a well-defined notion of a sequence of group elements $g_n$ going to infinity, indeed we have $\tau_{h_n} \to \infty$ if and only if $|h_n| \to \infty$. Given any set $K \subset H$, we let $GK := \{ gf: g \in G, f \in K \}$ be the orbit of $K$ under $G$.
\end{definition}

The $G$-invariance of the problem means in particular that the set $K_E$, being $G$-invariant, can no longer be compact (unless it consists only of $\{0\}$).  One might still hope that $K_E$ is an attractor in the sense of \eqref{tsim}, but this can be easily seen to be false (at least in the focusing NLS) by taking a stationary soliton solution $Q(x) e^{it}$ and applying a Galilean transform to create a travelling soliton which is not almost periodic.  The orbit of this travelling soliton is still almost periodic once one quotients out by the group $G$, so one might think to extend $K_E$ to cover solutions which are ``almost periodic modulo $G$'' (cf. \cite{compact}).  However, this is still not enough, as can be seen (at least heuristically) by considering \emph{multisoliton} solutions - the superposition of two or more diverging solitons.  See \cite{tsai-2} for details of how to construct such solutions forward-globally in time.  Observe that such solutions are not almost periodic even after quotienting out by $G$; in other words, there is no compact set $K \subset H$ such that the orbit $\{ u(t): t \in [0,+\infty) \}$ is contained in $GK$.  Thus a ``concentration compactness'' style definition of almost periodicity is needed, in order to account for the fact that solutions may be a superposition of components, each of which lives in a compact set after quotienting out by a different element of $G$.

To set this up properly requires some more notation.  

\begin{definition}[$G$-precompactness]\label{gpre}  If $K \subset H$ and $J \geq 0$ is an integer we let
$$J K := \{ f_1 + \ldots + f_J: f_1,\ldots,f_J \in K \}$$
denote the $J$-fold Minkowski sum of the set $K$, with the convention that $0K = \{0\}$.  We say that a set $E \subset H$ is \emph{$G$-precompact with $J$ components} if we have $E \subset J(G K)$ for some compact $K \subset H$ and $J \geq 1$.  We say that a solution $u: I \to H$ is \emph{$G$-almost periodic with $J$ components} if its orbit $\{ u(t): t \in I\}$ is $G$-precompact with $J$ components.
\end{definition}

\begin{example} Travelling soliton solutions are $G$-almost periodic with one component.  More generally, we expect multisoliton solutions formed by superimposing $J$ separated solitons to be $G$-almost periodic with $J$ components, provided that there is no radiation component whatsoever.  For other equivalent formulations of $G$-precompactness, see Proposition \ref{G-attractive}.
\end{example}

\begin{remark} The notion of being $G$-almost periodic with exactly one component is also known as being \emph{almost periodic modulo $G$}.  This type of almost periodicity is typically enjoyed by solitons, self-similar blowup solutions and by ``minimal blowup solutions''; see \cite{merlekenig}, \cite{compact} for further discussion.  Heuristically, we expect $G$-almost periodic solutions with $J$ components to be (nonlinear) superpositions of $J$ solutions which are almost periodic modulo $G$, but in the absence of an inverse scattering theory it is not entirely clear to the author what ``nonlinear superposition'' should mean.
\end{remark}

\begin{remark} The quantity $J = J(E)$ measures the maximum number of components associated to the asymptotic evolution of a solution of energy at most $E$, and thus we expect $J$ to grow at most linearly in $E$ (in light of Theorem \ref{scat}, we expect each non-radiating component to require a large amount of energy).  We will however not prove this claim here.
\end{remark}

We then have the following counterpart of Theorems \ref{attract-rad}, which we prove in Section \ref{nonradial-sec}.

\begin{theorem}[Non-radial compact attractor]\label{attract-nonrad} Let $E > 0$.  Then there exists a $G$-precompact closed NLS-invariant and $G$-invariant (i.e. translation-invariant) set $\K_E$ with $J = J(E) \geq 1$ components such that given any forward-global solution $u$ 
of energy at most $E$, there exists a unique \emph{radiation state} $u_+ \in H$ such that
\begin{equation}\label{tsim-nonrad}
\lim_{t \to +\infty} \dist_{H}( u(t) - e^{it\Delta} u_+, J \K_{E} ) = 0.
\end{equation}
In fact, we have a stronger statement: given any sequence $t_n$ of times going to $+\infty$, we have (after passing to a subsequence) a profile decomposition
\begin{equation}\label{profdecomp}
 u(t_n) = e^{it_n\Delta} u_+ + \sum_{j=1}^J \tau_{x_{j,n}} w_j + o_H(1)
 \end{equation}
where $w_j \in \K_E$, and $x_{j,n} \in \R^d$ obey the asymptotic separation condition
$$ \lim_{n \to \infty} |x_{j,n} - x_{j',n}| = \infty \hbox{ whenever } 1 \leq j < j \leq J.$$
\end{theorem}

\begin{remark} Informally, this theorem asserts that an arbitrary forward-global solution will asymptotically decouple into a radiation term $e^{it\Delta} u_+$, together with at most $J$ non-interacting (and widely separated) channels, each of which evolves within a $G$-precompact invariant set.  This latter set may still itself contain multiple components; this reflects the (rather strange) possibility that such a solution might have an infinite number\footnote{Imagine for instance a never-ending game of ``tennis'' in which one soliton is passed back and forth infinitely often between two other slowly diverging solitons.  While this type of scenario appears to be difficult to reconcile with conservation of momentum, it is not clear to the author how to rule it out completely.} of component interactions (such as collisions between components, fission of one component into multiple components, or fusion of multiple components into one), preventing an asymptotic resolution into non-interacting individual components.  One expects this to not be the case (at least for generic data) and so one might conjecture that one can take $\K_E$ to in fact be $G$-precompact with just one component; this would be consistent in particular with the soliton resolution conjecture.  Indeed, in view of Theorem \ref{attract-nonrad}, the soliton resolution conjecture is essentially equivalent to the assertion that $\K_E$ consists solely of soliton solutions, and is also equivalent to the assertion that the only $G$-almost periodic solutions are those which asymptotically resolve as the superposition of solitons.
\end{remark}

\begin{remark} Informally, \eqref{tsim-nonrad}, \eqref{profdecomp} are asserting that the non-radiating component of the solution $u(t)$ is localised to at most $J$ locations $x_1(t),\ldots,x_J(t)$.  In view of the soliton resolution conjecture, one expects these $x_j(t)$ to behave asymptotically linearly in $t$ for large $t$. There is some technical difficulty in formalising this statement rigorously, since at present one has the freedom at any time to perturb the points $x_j(t)$ by a displacement of $O(1)$, to permute the $x_j(t)$ with each other, and to replace two identical $x_j(t)$ with a single one (or vice versa), or to create and destroy dummy points (around which $u$ actually has very little mass).  In any event, the technology here does not seem sufficient to give such strong control on the $x_j(t)$, although it is likely that one can establish some sort of ``finite speed of propagation'' result which should allow one to take the $x_j(t)$ to be Lipschitz in $t$ (except when two trajectories coalesce, or one trajectory splits into two).
\end{remark}

We can now obtain non-radial extensions of Corollary \ref{petite-rad}:

\begin{corollary}[Petite conjecture, non-radial case I]\label{pet1} Let $u$ be a forward-global solution.  Then the following are equivalent:
\begin{itemize}
\item[(i)] $u$ is almost periodic.
\item[(ii)] $u$ is spatially localised near the origin in the sense that
$$\lim_{R \to +\infty} \limsup_{t \to +\infty} \int_{|x| > R} |u(t,x)|^2\ dx = 0.$$
\item[(iii)] $u$ is spatially localised near the origin in the sense that
$$ \lim_{R \to +\infty} \limsup_{t \to +\infty} \int_{|x| > R} |u(t,x)|^2 + |\nabla u(t,x)|^2\ dx = 0.$$
\end{itemize}
In either case we have $u_+ = 0$.
\end{corollary}

\begin{corollary}[Petite conjecture, non-radial case II]\label{pet2} Let $u$ be a forward-global solution.  Then the following are equivalent:
\begin{itemize}
\item[(i)] $u$ is $G$-almost periodic.
\item[(ii)] $u_+ = 0$ (i.e. $u$ is future non-radiating).
\item[(iii)] There exist functions $x_1, \ldots, x_J: \R^+ \to \R^d$ for some finite $J$ such that
$$ \lim_{R \to +\infty} \limsup_{t \to +\infty} \int_{|x-x_j(t)| > R \hbox{ for all } 1 \leq j \leq J} |u(t,x)|^2\ dx = 0.$$
\item[(iv)] There exist functions $x_1, \ldots, x_J: \R^+ \to \R^d$ for some finite $J$ such that
$$ \lim_{R \to +\infty} \limsup_{t \to +\infty} \int_{|x-x_j(t)| > R \hbox{ for all } 1 \leq j \leq J} |u(t,x)|^2 + |\nabla u(t,x)|^2\ dx = 0.$$
\end{itemize}
\end{corollary}

We establish these corollaries in Section \ref{nonradial-sec}.

\subsection{Organisation of the paper}

This paper is organised as follows.  After establishing some basic notation and estimates in Section \ref{notsec}, we discuss in Section \ref{motivation} a model instance of the ``double Duhamel trick'' which is absolutely essential to force compactness, and which relies heavily on the high dimension assumption $d \geq 5$.  We then in the remainder of the paper develop an increasingly sophisticated set of results concerning forward-global solutions.  Firstly, in Section \ref{loc-sec} we obtain some basic fixed-time and local-in-time estimates, basically arising from the energy bound and local (Strichartz) theory.  This, coupled with the decay of the fundamental solution, is already enough to construct the radiation state $u_+$ and the remaining bound state $v(t) := u(t) - e^{it\Delta} u_+$, which we do in Section \ref{compac-sec}.  Next, we exploit bilinear Strichartz estimates and the double Duhamel trick in Section \ref{freqloc-sec} to obtain localisation of $u(t)$ to frequencies of magnitude $\sim 1$.  Using this localisation, together with approximate finite speed of propagation and the double Duhamel trick again, we obtain also (in Section \ref{premspac-sec}) a preliminary localisation of $u(t)$ to boundedly many locations in physical space.  This preliminary localisation is already sufficient to extract a compact attractor in the radial case, which we do in Section \ref{radial-sec}.  In the non-radial case one needs to strengthen the localisation in a technical way using mass conservation, which do in Section \ref{fsl-sec}, before we can conclude the non-radial compactness results, which we establish in Section \ref{nonradial-sec}.

Finally, in Appendix \ref{dispersive} we collect some basic Strichartz-type estimate, while in Appendix \ref{appendix} we collect a number of basic facts about compact, precompact and $G$-precompact subsets of $H$ which we shall use several times in the main argument.

\subsection{Acknowledgements}

The author is supported by a grant from the Macarthur Foundation.  The author also thanks Tristan Roy for corrections and Igor Rodnianski for helpful comments, and in particular pointing out the connection to linear scattering theory.

\section{Notation}\label{notsec}

As mentioned in the introduction, we consider the parameters $d,p,\theta,C_0$ to be fixed, and allow all quantities to depend on these quantities.

We shall need four small exponents
$$ 1 \gg \eta_0 \gg \eta_1 \gg \eta_2 \gg \eta_3 > 0$$
where $\eta_0$ is assumed sufficiently small (depending on the above fixed parameters), $\eta_1$ is sufficiently small depending on $\eta_0$ (and the above fixed parameters), and so on down to $\eta_3$, which is sufficiently small depending on $\eta_0,\eta_1,\eta_2$ and the fixed parameters.

We also choose an arbitrary energy $E > 0$, which we now fix.

We use $X \lesssim Y$, $Y \gtrsim X$ or $X = O(Y)$ to denote the estimate $|X| \leq CY$, where $C$ depends only on the fixed parameters $d,p,\theta,C_0,F,E$ and on the exponents $\eta_0,\eta_1,\eta_2,\eta_3$.  If we need the constant $C$ to depend on other parameters too, we shall denote this by subscript, thus for instance $X \lesssim_\mu Y$, $Y \gtrsim_\mu$, or $X = O_\mu(Y)$ denotes the bound $|X| \leq C_\mu Y$ where $C$ can depend on $\mu$ as well as on the previous parameters.

We also use $O_H(Y)$ to denote an element of $H$ of norm $O(Y)$, and similarly with $H$ replaced by other normed vector spaces such as $L^2_x(\R^d)$.

We use $\langle x \rangle$ to denote the quantity $\langle x \rangle := (1 + |x|^2)^{1/2}$.  If $I$ is a time interval, we use $|I|$ to denote its length.

We shall need the following Littlewood-Paley projection operators.
Let $\varphi(\xi)$ be a bump function adapted to the ball $\{ \xi \in
\R^n: |\xi| \leq 2 \}$ which equals 1 on the ball $\{ \xi \in \R^n:
|\xi| \leq 1 \}$.  Define a \emph{dyadic number} to be any number $N \in 2^\Z$ of the form
$N = 2^j$ where $j \in \Z$ is an integer.  For each dyadic number $N$,
we define the Fourier multipliers
\begin{align*}
\widehat{P_{\leq N} f}(\xi) &:= \varphi(\xi/N) \hat f(\xi)\\
\widehat{P_{> N} f}(\xi) &:= (1 - \varphi(\xi/N)) \hat f(\xi)\\
\widehat{P_N f}(\xi) &:= (\varphi(\xi/N) - \varphi(2\xi/N)) \hat f(\xi).
\end{align*}
We similarly define $P_{<N}$ and $P_{\geq N}$.
We also define
$$ P_{M < \cdot \leq N} := P_{\leq N} - P_{\leq M} = \sum_{M < N' \leq N} P_{N'}$$
whenever $M < N$ are dyadic numbers.

We use the standard Sobolev spaces
$$ \| f\|_{W^{k,p}_x(\R^d)} := \sum_{j=0}^k \|\nabla^j f \|_{L^p_x(\R^d)}$$
for $1 \leq p \leq\infty$ and integers $k \geq 0$; throughout this paper $\nabla$ refers to the spatial gradient $\nabla_x$.

Given two norms $\|\|_V$ and $\|\|_W$, we write
$$ \|f\|_{V \cap W} := \|f\|_V + \|f\|_W.$$

\subsection{Special exponents}

Given any exponent $1 \leq q \leq \infty$ we let $q'$ be the dual exponent, defined by $\frac{1}{q}+\frac{1}{q'}=1$.

In the local theory we will need some special exponents $q_0,r_0,Q_0,Q,R$ which we now construct.

\begin{lemma}[Choice of exponents]\label{exponent}  
There exists an admissible pair $(q_0,r_0)$ with $q_0 > 2$,
exponents $2 < Q_0, Q < \frac{2d}{d-2}$, and an exponent $1 \leq R < \frac{2d}{d+4}$ such that
\begin{equation}\label{admis}
\frac{1}{r_0} + \frac{p-1}{Q_0} = \frac{1}{r'_0}
\end{equation}
and
\begin{equation}\label{irq}
\frac{1}{2} + \frac{p-1}{Q} = \frac{1}{R}.
\end{equation}
\end{lemma}

\begin{proof} This follows from the hypotheses $d \geq 5$, $1+\frac{4}{d} < p < 1+\frac{4}{d-2}$ and elementary algebra.  For instance, if $d=5$ and $p=2$ we can take $r_0=Q_0=3$, $R = 20/19$, $Q = 20/9$, and $q_0 = 12/5$; many other choices are of course possible.
\end{proof}

Henceforth we fix the exponents $q_0,r_0,Q_0,Q,R$ (which depend only on $d,p$)
and allow all quantities to depend on these exponents.

The exponents $q_0,r_0$ are useful for Strichartz iteration, see Lemma \ref{localst} below.
The significance of the $R$ exponent is the following:

\begin{lemma}[Fixed time estimate for $F(u)$]\label{ffix}  For any $u \in H$ we have
$$ \| F(u) \|_{W^{1,R}_x(\R^d)} \lesssim \|u\|_H^p.$$
\end{lemma}

\begin{proof}  
From \eqref{fpower}, \eqref{fpower-2}, \eqref{irq} and H\"older we have
$$ \| \nabla^j F(u(t)) \|_{L^R_x(\R^d)} \lesssim \| u(t) \|_{L^Q_x(\R^d)}^{p-1} \| \nabla^j u(t) \|_{L^2_x(\R^d)}$$
for $j=0,1$, and the claim follows from Sobolev embedding.
\end{proof}

\begin{remark}
Note that the hypothesis $d \geq 5$ was used in a crucial way in the argument, since the range $1 \leq R < \frac{2d}{d+4}$ is vacuous otherwise.  It is very important for us that $R$ lie in this range, as it implies that the relevant dispersive inequality \eqref{energy-dispersive} enjoys a decay which is strictly better than $|t|^{-2}$.
\end{remark}

\section{Motivation: using dispersive estimates to exclude resonances}\label{motivation}

In the spectral theory of \emph{linear} Schr\"odinger operators $-\Delta + V$ with potential, it is well known that in high dimension $d \geq 5$ and with rapidly decaying $V$ that there are no resonances, i.e. weakly growing non-$L^2(\R^d)$ solutions $\phi$ to the eigenfunction equation $(-\Delta + V) \phi = E \phi$.  This result is usually established by exploiting the strong \emph{spatial} decay (better than $|x|^{-d/2}$) of the resolvent kernels $(-\Delta-E \pm i\eps)^{-1}$.  In this section we briefly (and informally) indicate how, via a Fourier transform in time and a ``double Duhamel trick'', one can establish the same result instead using the strong \emph{time} decay (better than $t^{-2}$) of the Schr\"odinger kernels of $e^{it\Delta}$.  We give this alternate derivation here because the same double Duhamel trick shall be relied upon heavily throughout the remainder of the paper; furthermore,
this exclusion of resonances result for linear Schr\"odinger equations can be viewed as somewhat analogous to the compact attractor results for nonlinear Schr\"odinger equations presented here.  It seems a worthwhile task to try to further develop the analogy between linear spectral and scattering theory and nonlinear scattering theory.

As this section is only included for motivational purposes (and will not be directly used elsewhere in this paper) 
we shall work non-rigorously and imprecisely.  Let $f \in L^1(\R^d) \cap L^2(\R^d)$.  For any real energy $E$, we \emph{formally} have the identities
$$ (-\Delta - E)^{-1} f = -i \int_0^\infty e^{iEt} e^{it\Delta} f \ dt$$
and
$$ (-\Delta - E)^{-1} f = i \int_{-\infty}^0 e^{iEt'} e^{it'\Delta} f \ dt';$$
to make these identities well-defined and rigorous one needs the limiting absorption principle, which we will not discuss here.  Taking inner products we obtain the formal ``double Duhamel'' identity
$$ \| (-\Delta - E)^{-1} f \|_{L^2(\R^d)}^2
= - \int_0^\infty \int_{-\infty}^0 e^{iE(t-t')} \langle e^{it\Delta} f, e^{it'\Delta} f \rangle\ dt dt'$$
and hence by the triangle inequality
$$ \| (-\Delta - E)^{-1} f \|_{L^2(\R^d)}^2
\leq \int_0^\infty \int_{-\infty}^0 |\langle e^{it\Delta} f, e^{it'\Delta} f \rangle|\ dt dt'.$$
Now we have two estimates for the inner product.  On the one hand, by Cauchy-Schwarz we can bound this inner product by
$\|f\|_{L^2(\R^d)}^2$.  On the other hand, by the dispersive inequality we can bound the inner product by $O( \frac{1}{|t-t'|^{d/2}} \|f\|_{L^1(\R^d)}^2 )$.  Putting the two bounds together we obtain
$$ \| (-\Delta - E)^{-1} f \|_{L^2(\R^d)}^2
\lesssim \int_0^\infty \int_{-\infty}^0 \frac{1}{\langle t-t'\rangle^{d/2}} \|f\|_{L^1(\R^d) \cap L^2(\R^d)}^2\ dt dt'.$$
But when $d \geq 5$, the integrand is convergent, and we conclude that
$$ \| (-\Delta - E)^{-1} f \|_{L^2(\R^d)} \lesssim \|f\|_{L^1(\R^d) \cap L^2(\R^d)}.$$
This estimate can be used to exclude resonances; indeed, if $(-\Delta + V)\phi = E \phi$, then
we formally have $\phi = (-\Delta-E)^{-1} (-V\phi)$.  If $\phi$ is slowly growing and $V$ is sufficiently rapidly decreasing, then $V\phi \in L^1(\R^d) \cap L^2(\R^d)$, and we conclude that $\phi \in L^2(\R^d)$, thus $\phi$ is a genuine eigenfunction rather than a resonance.

\section{Local estimates}\label{loc-sec}

Let $u$ be a forward-global solution of energy at most $E > 0$.  In this section we allow all implied
constants to depend on $E$.

In this section we obtain estimates on $u$ which are either fixed-time or local in time.

\subsection{Fixed-time estimates}

From \eqref{energy-bound} and our convention to suppress dependence on $E$ we have
\begin{equation}\label{energy}
\| u(t) \|_{H} \lesssim 1 \hbox{ for all } t \in [0,+\infty)
\end{equation}
and hence by Sobolev embedding we have
\begin{equation}\label{sob}
\| u(t) \|_{L^q_x(\R^d)} \lesssim 1 \hbox{ for all } 2 \leq q \leq \frac{2d}{d-2} \hbox{ and } t \in [0,+\infty).
\end{equation}
From Lemma \ref{ffix} we also conclude
\begin{equation}\label{ffix-est}
\| F(u(t)) \|_{W^{1,R}_x(\R^d)} \lesssim 1 \hbox{ for all } t \in [0,+\infty).
\end{equation}

\subsection{Local-in-time estimates}

Next we establish a (very standard) local-in-time Strichartz estimate.  

\begin{lemma}[Local Strichartz control]\label{localst}  For any time interval $I \subset [0,+\infty)$ and any admissible pair of exponents $(q,r)$ we have
\begin{equation}\label{nqr}
 \| u \|_{L^q_t W^{1,r}_x( I \times \R^d ) } \lesssim \langle |I| \rangle^{1/q}
 \end{equation}
as well as the nonlinearity estimate
\begin{equation}\label{naff} 
\| F(u) \|_{L^{q'_0}_t W^{1,r'_0}_x( I \times \R^d ) } \lesssim \langle |I| \rangle^{1/q'_0}.
\end{equation}
(Recall $q_0,r_0$ were fixed by Lemma \ref{exponent}.)
\end{lemma}

\begin{proof}  As the arguments here are very standard (see e.g. \cite{caz}, \cite{cwI}, \cite{katounique}) we give only a sketch here.
By subdividing $I$ it suffices to prove this claim in the case when the interval $I$ has length much smaller than $1$.  Let us first work formally, assuming \emph{a priori}
that all norms appearing below are finite.  Define
the quantity 
$$ X := \| u \|_{L^{q_0}_t W^{1,r_0}_x(I \times \R^d)}.$$
From \eqref{homog}, \eqref{retarded}, \eqref{energy} we have
$$ X \lesssim 1 + \| F(u) \|_{L^{q'_0}_t W^{1,r'_0}_x(I \times \R^d)}.$$
From \eqref{fpower}, \eqref{fpower-2} we have $\nabla^j F(u) = O( |u|^{p-1} |\nabla^j u| )$ for $j=0,1$.
From H\"older's inequality, together with \eqref{admis}, we thus see that
\begin{equation}\label{sumj}
\| F(u) \|_{L^{q'_0}_t W^{1,r'_0}_x(I \times \R^d)}
\lesssim |I|^{1/q_0 - 1/q'_0} \| u \|_{L^\infty_t L^{Q_0}_x(I \times \R^d)}^{p-1} X.
\end{equation}
Applying \eqref{sob} we conclude the \emph{a priori} estimate
$$ X \lesssim 1 +|I|^{1/q_0 - 1/q'_0} X.$$
Since $q_0 < 2$, we thus obtain a bound of the form
$$ X \leq O(1) + \frac{1}{2} X$$
if we make $|I|$ sufficiently small.  This gives a bound $X = O(1)$ provided that $X$ is finite.  If we then set up a standard Picard iteration scheme to construct the solution $u$ on $I \times \R^d$ (starting from initial condition $u(t_0)$ for some $t_0 \in I$), exploiting the uniqueness theory in \cite{katounique}, and adapt the above argument to the iterates, we conclude that we indeed do have the bound $X = O(1)$.  Applying the Strichartz estimates (and \eqref{sumj}) one last time gives the claim.
\end{proof}

We have a bilinear refinement of the above estimate which shows that interactions between widely different frequencies is weak.

\begin{corollary}[Bilinear estimate]\label{bil}  For any time interval $I \subset [0,+\infty)$, any $0 \leq \delta < 1/2$, and any $N, M > 0$ we have
$$ \| |u_N| |u_M| \|_{L^2_t L^2_x(I \times \R^d)} \lesssim_\delta \langle |I| \rangle^{1/2}
\frac{M^\delta}{N^\delta} \frac{M^{\frac{d-2}{2}}}{\langle N \rangle \langle M \rangle}.$$ 
\end{corollary}

Of course, this corollary is most effective in the regime $N \geq M$; in the opposite regime $N \leq M$, one should swap the roles of $N$ and $M$.

\begin{proof} Again, by subdividing time we may assume $|I|$ to be less than $1$.  One now applies
Theorem \ref{bilst} for some $t_0 \in I$ followed by Lemma \ref{localst} (and \eqref{energy}) and Plancherel's theorem, observing that $(i\partial_t + \Delta) u_N = P_N F(u)$ and similarly for $P_M u$.
\end{proof}

Corollary \ref{bil} gives some additional regularity\footnote{This ``smoothing effect'' in the nonlinearity in the energy-subcritical case has also been exploited in the global well-posedness theory below the energy norm, see \cite{borg:book}.  The proposition below is also very closely related to some estimates in \cite{ckstt:8}.} on $F(u)$ beyond the one derivative that appears
in Lemma \ref{localst}:

\begin{proposition}[Smoothing effect]\label{smooth-effect}  There exists $\sigma > 0$ (depending only on $d,p,q_0,r_0$) such that we have the smoothing estimate
\begin{equation}\label{smooth}
 \| P_N F(u) \|_{L^{q'_0}_t L^{r'_0}_x( I \times \R^d ) } \lesssim \langle N \rangle^{-1-\eta_1} \langle |I| \rangle^{1/q'_0}.
 \end{equation}
\end{proposition}

\begin{proof} As before we may take $|I| \leq 1$.  We may also take $N \ge 1$ as the claim follows from Lemma \ref{localst} otherwise.  For brevity we shall omit the domain $I \times \R^d$ in all the norms in this proof.  Applying a derivative, our task is to show that
\begin{equation}\label{smooth-2}
 \| P_N (F'(u) \nabla u) \|_{L^{q'_0}_t L^{r'_0}_x} \lesssim N^{-\eta_1}
\end{equation}
where we think of $F'(u)$ as a $2 \times 2$ matrix and $\nabla u$ as a $2 \times 1$ column vector.

We shall use a rather crude paradifferential method\footnote{We will not need the full power of the paradifferential calculus here because we are not seeking an optimal value of $\eta_1$.  The reader may wish to work out the exponents explicitly in the model case $(d,p) = (5,2)$, which allows for some slight simplifications.}.  We split $u = u_\lo + u_\hi$, where
$u_\lo := u_{<N^{\eta_0}}$ and $u_\hi := u_{\geq N^{\eta_0}}$.  From \eqref{fpower-3} we have
$$ F'(u) = F'(u_\lo) + O( |u_\hi|^\theta (|u_\lo| + |u_\hi|)^{p-1-\theta} ).$$
Consider the contribution of the error term to \eqref{smooth-2}.  We discard the $P_N$ and use H\"older (and \eqref{admis}) followed by \eqref{nqr}, \eqref{sob} to estimate this contribution by
$$ \lesssim \| u_\hi \|_{L^\infty_t L^{Q_0}_x}^{\theta} \| |u_\lo| + |u_\hi| \|_{L^\infty_t L^{Q_0}_x}^{p-1-\theta} \| \nabla u \|_{L^{q_0}_t L^{r_0}_t}
\lesssim \| u_\hi \|_{L^\infty_t L^{Q_0}_x}^{\theta}.$$
But by construction, $Q_0$ is strictly less than the endpoint Sobolev exponent $\frac{2d}{d-2}$, so from \eqref{energy} and Bernstein's inequality we have
$$ \| u_\hi \|_{L^\infty_t L^{Q_0}_x} \lesssim (N^{\eta_0})^{-c}
$$
for some $c$ depending only on $Q_0$ and $d$.  This gives a net contribution of $N^{-c \theta \eta_0}$ to \eqref{smooth-2}, which is acceptable.

Having disposed of the error term, it remains to show that
$$  \| P_N (F'(u_\lo) \nabla u) \|_{L^{q'_0}_t L^{r'_0}_x} \lesssim N^{-\eta_1}.$$
We split $\nabla u = \nabla u_{<N/100} + \nabla u_{\geq N/100}$, where $u_{<N/100} = P_{<N/100} u$ and $u_{\geq N/100} = P_{\geq N/100} u$.  We consider first the contribution of $\nabla u_{<N/100}$.  Then we may freely replace $F'(u_\lo)$ by $P_{>N/10} F'(u_\lo)$.  By H\"older and \eqref{admis} followed by \eqref{nqr} we may thus estimate this contribution by
\begin{align*}
\| P_N (F'(u_\lo) \nabla u_{<N/100}) \|_{L^{q'_0}_t L^{r'_0}_x} 
&=
\| P_N ((P_{>N/10} F'(u_\lo)) \nabla u_{<N/100}) \|_{L^{q'_0}_t L^{r'_0}_x} \\
&
\lesssim \| P_{>N/10} F'(u_\lo) \|_{L^\infty_t L^{(p-1)Q_0}_x} \| \nabla u_{\geq N/100} \|_{L^{q_0}_t L^{r_0}_x} \\
&\lesssim \| P_{>N/10} F'(u_\lo) \|_{L^\infty_t L^{(p-1)Q_0}_x}.
\end{align*}
From the rapid decrease and cancellation properties of the convolution kernel of $P_{>N/10}$ and Minkowski's inequality (and the triangle inequality) we have the H\"older regularity-type bound
$$ \| P_{>N/10} F'(u_\lo) \|_{L^\infty_t L^{(p-1)Q_0}_x} \lesssim 
\sup_{|h| \leq 1/N} \| \tau_h F'(u_\lo) - F'(u_\lo) \|_{L^\infty_t L^{(p-1)Q_0}_x}.$$
But from \eqref{fpower-3} we have the pointwise estimate
$$|\tau_h F'(u_\lo) - F'(u_\lo)| \lesssim |\tau_h u_\lo - u_\lo|^\theta (|\tau_h u_\lo| + |u_\lo|)^{p-1-\theta}.$$
Combining all this with H\"older's inequality followed by \eqref{sob} we conclude that 
\begin{align*}
\| P_N (F'(u_\lo) \nabla u_{<N/100}) \|_{L^{q'_0}_t L^{r'_0}_x} 
&\lesssim \sup_{|h| \leq 1/N} \| \tau_h u_\lo - u_\lo \|_{L^\infty_t L^{Q_0}_x}^\theta
\| |\tau_h u_\lo| + |u_\lo| \|_{L^\infty_t L^{Q_0}_x}^{p-1-\theta}\\
&\lesssim \sup_{|h| \leq 1/N} \| \tau_h u_\lo - u_\lo \|_{L^\infty_t L^{Q_0}_x}^\theta
\end{align*}
By the fundamental theorem of calculus and Minkowski's inequality we conclude
$$ \| P_N (F'(u_\lo) \nabla u_{<N/100}) \|_{L^{q'_0}_t L^{r'_0}_x} \lesssim \left(\frac{1}{N} \| \nabla u_\lo \|_{L^\infty_t L^{Q_0}_x}\right)^\theta.$$
But by \eqref{sob} and Littlewood-Paley theory the right-hand side is$O(\frac{1}{N} N^{\eta_0})^{\theta}$, which is  acceptable.

It remains to show that
$$  \| P_N (F'(u_\lo) \nabla u_{\geq N/100}) \|_{L^{q'_0}_t L^{r'_0}_x} \lesssim N^{-\eta_1}.$$
We discard the $P_N$ and apply \eqref{fpower-2} to estimate the left-hand side by
$$ \| |u_\lo|^{p-1} |\nabla u_{\geq N/100}| \|_{L^{q'_0}_t L^{r'_0}_x}.$$
Now observe that from H\"older's inequality and \eqref{sob}, \eqref{nqr} (as in \eqref{sumj})
we have
$$ \| |u_\lo|^{p-1} |\nabla u_{\geq N/100}| \|_{L^{q'_0}_t L^{r'_0}_x}
\lesssim \| u_\lo \|_{L^\infty_t L^{Q_0}_x}^{p-1} \| \nabla u_{\geq N/100} \|_{L^{q_0}_t L^{r_0}_x} \lesssim 1.$$
Indeed, because there is some ``room'' in the use of non-endpoint Sobolev embedding and a H\"older inequality in time, we in fact see that we have the more general estimates
$$ \| |u_\lo|^{\tilde p-1} |\nabla u_{\geq N/100}| \|_{L^{\tilde q'_0}_t L^{\tilde r'_0}_x} \lesssim 1$$
for all $(\tilde p, \tilde q_0, \tilde r_0)$ in a sufficiently small neighbourhood of $(p, q_0, r_0)$.
On the other hand, from a dyadic decomposition argument followed by Corollary \ref{bil} (swapping the roles of $N$ and $M$) 
yields the bound
$$ \| |u_\lo| |\nabla u_{\geq N/100}| \|_{L^2_t L^2_x} \lesssim N^{-0.1}$$
(say) since $\eta_0$ is small.  The claim then follows by a suitable application of H\"older's inequality (or interpolation).
\end{proof}

\begin{remark}
The above smoothing effect will be crucial in obtaining a high frequency decay estimate which will form one component of the compactness result.  (The other components are a low frequency decay estimate and a spatial localisation estimate.)
\end{remark}

\section{Preliminary asymptotic analysis}\label{compac-sec}

We now combine the above local estimates with the strong decay of the fundamental solution in high dimension to control the asymptotic frequency distribution.  We begin with

\begin{lemma}[Weak scattering]\label{weaks}  Let $u$ be a forward-global solution of energy at most $E$.
The elements $e^{-it\Delta} u(t)$ of $H$ are weakly convergent as $t \to +\infty$.  In other words, for any function $\phi \in H$, the sequence $\langle e^{-it\Delta} u(t), \phi \rangle_{H}$ is convergent as $t \to +\infty$.
\end{lemma}

\begin{proof} By \eqref{energy-bound} we see that the $e^{-it\Delta} u(t)$ are uniformly bounded in $H$.  Thus it will suffice by the usual limiting arguments to verify the weak convergence in the case where $\phi \in C^\infty_0(\R^d)$ is a test function.  It suffices to show that
$$ \lim_{t_1, t_2 \to +\infty} \langle e^{-it_1 \Delta} u(t_1) - e^{-it_2 \Delta} u(t_2), \phi \rangle_{H} = 0.$$ 
But from Duhamel's formula we have
$$ e^{-it_1 \Delta} u(t_1) - e^{-it_2 \Delta} u(t_2) = - \int^{t_1}_{t_2} e^{-it\Delta} F(u(t))\ dt$$
and so it will suffice by Minkowski's inequality to show that
$$ \lim_{t_1, t_2 \to +\infty} \int_{t_1}^{t_2} |\langle F(u(t)), e^{it\Delta} \phi \rangle_{H}|\ dt = 0.$$
By Lemma \ref{ffix} and H\"older's inequality, it suffices to show that
$$ \limsup_{t_1, t_2 \to +\infty} \int_{t_1}^{t_2} \| e^{it\Delta} \phi \|_{W^{1,R'}_x(\R^d)}\ dt = 0.$$
But from \eqref{energy-dispersive} and the hypothesis that $\phi$ is a test function we have
$$ \| e^{it\Delta} \phi \|_{W^{1,R'}_x(\R^d)} \lesssim_\phi |t|^{d(\frac{1}{R} - \frac{1}{2})}.$$
But since $R < \frac{2d}{d+4}$, the exponent here is greater than $2$, and the claim certainly follows.
\end{proof}

We conclude

\begin{proposition}[Preliminary decomposition]\label{pre}  Let $u$ be a forward-global solution with energy $E > 0$.  Then there
exists a unique decomposition
\begin{equation}\label{up}
 u(t) = e^{it\Delta} u_+ + v(t)
 \end{equation}
where $u_+ \in H$ with 
\begin{equation}\label{uplus-bound}
\|u_+\|_H^2 \leq E
\end{equation}
and
\begin{equation}\label{vbound}
\|v(t)\|_H \leq 2 \sqrt{E}
\end{equation}
for all $t \in [0,+\infty)$.  Furthermore, 
$v$ is \emph{weakly bound} in the sense that
\begin{equation}\label{vweak}
\wlim e^{-it\Delta} v(t) = 0
\end{equation}
or equivalently that
\begin{equation}\label{uplus}
u_+ = \wlim e^{-it\Delta} u(t),
\end{equation}
where the limit is in the weak topology.
We also have the Duhamel identities
\begin{equation}\label{uforward}
v(t) = e^{it\Delta} [u(0)-u_+] - i \int_0^t e^{i(t-t'')\Delta} F(u(t''))\ dt''
\end{equation}
and
\begin{equation}\label{ubackward}
 v(t) = i \Wlim \int_t^T e^{i(t-t')\Delta} F(u(t'))\ dt'.
\end{equation}
\end{proposition}

\begin{proof}  We define $u_+$ by \eqref{uplus}, which exists by Lemma \ref{weaks}.
From \eqref{energy} and the weak closure of the unit ball of $H$ we have the desired $H$-norm
bound on $u_+$, and the bound on $v$ follows by the triangle inequality.
We then write
$$ v(t) := u(t) - e^{it\Delta} u_+,$$
from which \eqref{vweak} and \eqref{uforward} are automatic.
On the other hand, from Duhamel's formula (backward in time) we have
$$ u(t) = e^{it\Delta} e^{-iT\Delta} u(T) + i \int_t^T e^{i(t-t')\Delta} F(u(t'))\ dt'$$
for any $T > t$; we obtain \eqref{ubackward}.

For the uniqueness, observe that on applying $e^{-it\Delta}$ to \eqref{up} and taking weak limits (using \eqref{vweak}) we conclude \eqref{uplus}, which defines $u_+$ and hence $v(t)$ uniquely.
\end{proof}

We shall refer to $v$ as the \emph{weakly bound component} of $u$.  In the next few sections we analyse this component further.  A key strategy will be to take inner products between \eqref{uforward} and \eqref{ubackward} in order to control $v$ in $L^2$ by a double Duhamel integral.  It is here that the high dimension hypothesis $d \geq 5$ becomes decisive, basically because the integral $\int_0^t \int_t^\infty \frac{1}{\langle t'-t''\rangle^{d/2}}\ dt' dt''$ is bounded uniformly in $t$ in that case (cf. Section \ref{motivation}).  

\begin{remark} The above analysis only requires the fundamental solution to decay faster than $1/t$.  As such, it works in three and higher dimensions, not just in five and higher dimensions.  For instance, the $d=3,p=3$ case of the above results were established in \cite{tao:compact}.
\end{remark}

\begin{remark} In practice, the radiation term $e^{it\Delta} u_+$ (as well as the linear solution $e^{it\Delta} u_0$) will be asymptotically negligible as $t \to +\infty$, thanks to Lemma \ref{rl}.
\end{remark}

\begin{remark} From \eqref{vweak} we have the asymptotic orthogonality estimate
$$ \lim_{t \to +\infty} \| v(t) + e^{it\Delta} u_+ \|_H^2 - \|v(t)\|_H^2 - \| e^{it\Delta} u_+\|_H^2 = 0$$
and hence we have an asymptotic decoupling of energy
$$ \|u_+\|_H^2 + \limsup_{t \to +\infty} \|v(t)\|_H^2 = \limsup_{t \to +\infty} \|u(t)\|_H^2 \leq E.$$
There is a similar decoupling for the mass $M()$ and Hamiltonian $H()$ (cf. \cite{tao:compact}).  We will not use these facts here.
\end{remark} 

It will be useful to note that the weakly bound component $v$ itself is an approximate solution to NLS:

\begin{lemma}[Weakly bound states are approximate solutions]\label{approx}  Let $T \in \R$, let $u$ be a forward-global solution, and let $v$ be the weakly bound state.  Then for all sufficiently late times $t_0$ (depending on $T,u$) we have
$$ S(T) v(t_0) = v(T+t_0) + o_H(1)$$
where $o_H(1)$ goes to zero in $H$ norm as $t_0 \to +\infty$.
\end{lemma}

\begin{proof}  
Fix $T$, and let $t_0$ be a sufficiently late time.  Let $I$ be the interval $I := [t_0,t_0+T]$.
Observe that $v$ solves the forced NLS
$$ i v_t + \Delta v = F(v) + [F(v+e^{it\Delta} u_+) - F(v)]$$
on $I$.  

From \eqref{uplus-bound} and Lemma \ref{homog} we have
\begin{equation}\label{sujo}
 \| e^{it\Delta} u_+ \|_{L^{q_0}_t W^{1,r_0}_x( I \times \R^d ) } = o(1)
\end{equation}
where $o(1)$ is a quantity which goes to zero as $t_0 \to \infty$.  Similarly, from Lemma \ref{rl} we have
\begin{equation}\label{sujo-2}
 \| e^{it\Delta} u_+ \|_{L^\infty_t L^Q_x( I \times \R^d ) } = o(1)
\end{equation}
while from \eqref{uplus-bound} we also have
$$ \| e^{it\Delta} u_+ \|_{C^0_t H^1_x(I \times \R^d)} \lesssim 1.$$
From Lemma \ref{localst} and the triangle inequality we conclude
\begin{equation}\label{sumjv}
\| v \|_{L^{q_0}_t W^{1,r_0}_x \cap C^0_t H^1_x(I \times \R^d)} \lesssim_T 1.
\end{equation}
Now let us compute the quantity
$$ X := \|F(v+e^{it\Delta} u_+) - F(v) \|_{L^{q'_0}_t W^{1,r'_0}_x(I \times \R^d)}.$$
From \eqref{fpower}, \eqref{fpower-2}, \eqref{fpower-3} we have the pointwise bound
\begin{align*}
|\nabla^j (F(v+e^{it\Delta} u_+) - F(v))| &\lesssim
 (|v| + |e^{it\Delta} u_+|)^{p-1} |\nabla^j e^{it\Delta} u_+|\\
&\quad + (|v| + |e^{it\Delta} u_+|)^{p-1-\theta} |e^{it\Delta} u_+|^{\theta} |\nabla^j v|
\end{align*}
for $j=0,1$, and so by H\"older and \eqref{admis}, \eqref{sumjv}, \eqref{sujo}, \eqref{sujo-2} we have $X=o(1)$.
 The claim now follows (for sufficiently late times $t_0$) by invoking Lemma \ref{perturb}.
\end{proof}

\section{Frequency localisation}\label{freqloc-sec}

We now localise the weakly bound component $v$ of the forward-global solution $u$ in frequency.

\begin{proposition}[Asymptotic frequency localisation of energy]\label{asymloc}  Let $u$ be a forward-global solution of energy at most $E$, and let $v$ be the weakly bound component of $u$.
Then we have
\begin{equation}\label{pnlow}
\limsup_{t \to +\infty} \sup_{N \leq 1} N^{-\eta_3} \| P_{\leq N} v(t) \|_{H} \lesssim 1
\end{equation}
and
\begin{equation}\label{pnbig}
\limsup_{t \to +\infty} \sup_{N \geq 1} N^{\eta_3} \| P_{\geq N} v(t) \|_{H} \lesssim 1.
\end{equation}
\end{proposition}

\begin{remark} From Proposition \ref{attractive} in the appendix we see that obtaining these types of frequency localisations are an important step towards compact attractor results such as Theorem \ref{attract-rad} and Theorem \ref{attract-nonrad}, although one of course will also need spatial localisation results to conclude these theorems.  It is significantly easier to establish localisation in frequency rather than in space, largely due to the fact that the linear propagators $e^{it\Delta}$ preserve the former but not the latter.
\end{remark}

\begin{proof}  
The main idea will be to play off the Duhamel formulae \eqref{ubackward} and \eqref{uforward} against each other (cf. Section \ref{motivation}). 

Fix $\eps > 0$.  Since test functions are dense in $H$, we can find
a test function $u_\eps \in C^\infty_0(\R^d)$ such that 
$$ u(0) - u_+ = u_\eps + O_{H}(\eps^2)$$
and hence from \eqref{uforward}
\begin{equation}\label{uf}
v(t) = e^{it\Delta} u_\eps - i \int_0^t e^{i(t-t'')\Delta} F(u(t''))\ dt'' + O_{H}(\eps^2).
\end{equation}

Now we prove \eqref{pnlow}.  It suffices to show that
\begin{equation}\label{pnlow-a}
\| P_{\leq N} v(t) \|_{H} \lesssim N^{\eta_3} +\eps
\end{equation}
for all sufficiently late times $t$ and $N \geq 1$, where ``sufficiently late'' can depend on $E,u,\eps$ and the fixed parameters but is uniform in $N$.

Applying $P_{\leq N}$ to \eqref{uf} and to \eqref{ubackward} we obtain
$$
P_{\leq N} v(t) = i \Wlim \int_t^T e^{i(t-t')\Delta} P_N F(u(t'))\ dt' 
$$
and
$$
P_{\leq N} v(t) = P_{\leq N} e^{it\Delta} u_\eps - i \int_0^t e^{i(t-t'')\Delta} P_N F(u(t'')\ dt'' + O_{H}(\eps^2).
$$
On the other hand, from \eqref{vbound} we have
$$ P_{\leq N} v(t) = O_{H}(1)$$
and hence
\begin{align*}
\| P_{\leq N} v(t) \|_{H}^2 
&= |\langle P_{\leq N} v(t), P_{\leq N} v(t) \rangle_{H}| \\
&\leq |\langle P_{\leq N} v(t), P_{\leq N} e^{it\Delta} u_\eps - i \int_0^t e^{i(t-t'')\Delta} P_N F(u(t'')\ dt'' 
\rangle_{H}| + O(\eps^2) \\
&= \langle i \Wlim \int_t^T e^{i(t-t')\Delta} P_{\leq N} F(u(t'))\ dt' , \\
&\quad\quad P_{\leq N} e^{it\Delta} u_\eps - i \int_0^t e^{i(t-t'')\Delta} P_{\leq N} F(u(t'')\ dt'' 
\rangle_{H} + O(\eps^2) \\
&\leq \int_t^T |\langle e^{i(t-t')\Delta} P_{\leq N} F(u(t')), P_{\leq N} e^{it\Delta} u_\eps \rangle_{H}|\ dt' \\
&\quad + |\int_t^T \int_0^t Y_N(t',t'')\ dt'' dt'| + O(\eps^2)
\end{align*}
for sufficiently large $T$ (depending on all other quantities)
where $Y_N(t',t'')$ is the quantity
$$ Y_N(t',t'') := \langle e^{i(t-t')\Delta} P_{\leq N} F(u(t')), e^{i(t-t'')\Delta} P_{\le N} F(u(t'')) \rangle_{H}.$$
For the first integral, we use integration by parts and H\"older's inequality, followed by \eqref{ffix-est}, to write
$$ |\langle e^{i(t-t')\Delta} P_N F(u(t')), P_N e^{it\Delta} u_\eps \rangle_{H}|
\lesssim \| P_N^2 e^{it' \Delta} u_\eps \|_{W^{2,R'}_x}.$$
From \eqref{energy-dispersive} and Sobolev embedding we have
$$ \| P_N^2 e^{it' \Delta} u_\eps \|_{W^{2,R'}_x} \lesssim_{u_\eps} \frac{1}{\langle t' \rangle^{d (\frac{1}{R}-\frac{1}{2})}}.$$
Since $R < \frac{2d}{d+2}$, the exponent of $\langle t' \rangle$ here is more than $2$, and we obtain
$$ \int_t^\infty |\langle e^{i(t-t')\Delta} P_N F(u(t')), P_N e^{it\Delta} u_\eps \rangle_{H}| = O(\eps^2)$$
if $t$ is a sufficiently late times depending on $E$, $\eps$, and $u_\eps$ (but \emph{not} depending on $N$).
We conclude that
$$ \| P_{\leq N} v(t) \|_{H}^2 \lesssim \eps^2 + |X_{\leq N}|.$$

We note for future reference that the only properties of the linear operator $P_{\leq N}$ that we have used so far are that it commutes with Fourier multipliers and is bounded on all Lebesgue and Sobolev spaces.  In particular, the arguments apply with $P_{\leq N}$ replaced by any other Littlewood-Paley multiplier.

To conclude the proof of \eqref{pnlow-a} and hence \eqref{pnlow} it suffices to show that
$$ \int_t^\infty \int_0^t |Y_N(t',t'')|\ dt'' dt' \lesssim N^{-\eta_3}.$$
We use the low frequency hypothesis $N \leq 1$ to estimate
$$ |Y(t',t'')| \lesssim \left|\langle P_{\leq N} F(u(t')), e^{i(t'-t'')\Delta} P_{\leq N} F(u(t'')) \rangle_{L^2_x(\R^d)}\right|.
$$
From \eqref{energy-dispersive} followed by \eqref{ffix-est}  we thus have
\begin{align*}
|Y(t',t'')|
&\lesssim \frac{\| P_{\leq N} F(u(t')) \|_{L^R_x(\R^d)} \| P_{\leq N} F(u(t'')) \|_{L^R_x(\R^d)}}{(t'-t'')^{d(\frac{1}{R}-\frac{1}{2})}}  \\
&\lesssim \frac{1}{(t'-t'')^{d(\frac{1}{R} - \frac{1}{2})}}.
\end{align*}
On the other hand, from \eqref{energy-dispersive}, \eqref{ffix-est} and Bernstein's inequality we have
\begin{align*}
|Y(t',t'')| &\lesssim \| P_{\leq N} F(u(t'))\|_{L^2_x(\R^d)} \| P_{\leq N} F(u(t''))\|_{L^2_x(\R^d)} \\
&\lesssim N^{d(\frac{1}{R}-\frac{1}{2})} \| F(u(t'))\|_{L^R_x(\R^d)} N^{d(\frac{1}{R}-\frac{1}{2})} \| F(u(t''))\|_{L^R_x(\R^d)} \\
&\lesssim N^{d(\frac{2}{R}-1)}.
\end{align*}
Putting all this together we obtain
$$ \int_t^\infty \int_0^t |Y_N(t',t'')|\ dt'' dt' \lesssim \int_t^\infty \int_0^t \min\left(\frac{1}{(t'-t'')}, N^2\right)^{d(\frac{1}{R}-\frac{1}{2})}\ dt'' dt'.$$
Since $R < \frac{2d}{d+2}$, the exponent $d(\frac{1}{R}-\frac{1}{2})$ is strictly greater than two\footnote{This is the one place where we crucially rely on the hypothesis that the dimension is at least five, as our arguments very much need a dispersive estimate which decays faster than $1/t^2$ in order to obtain localisation in the double Duhamel integral.}, and the claim \eqref{pnlow-a} follows (for $\delta$ sufficiently small).

Now we establish \eqref{pnbig}.  We argue as before and conclude that for sufficiently late times $t$ (uniformly in $N$) as before we have
$$ \| P_{\geq N} v(t) \|_{H}^2 \lesssim \eps^2 + |\int_t^T \int_0^t Z_N(t',t'')\ dt'' dt'| 
$$
for sufficiently large $T$,
where 
$$ Z_N(t',t'') := \langle e^{i(t-t')\Delta} P_{\geq N} F(u(t')), e^{i(t-t'')\Delta} P_{\geq N} F(u(t'')) \rangle_{H}.$$
It thus suffices to show that
$$ |\int_t^T \int_0^t Z_N(t',t'')\ dt'' dt'| \lesssim N^{-\eta_3}.$$
Let us first dispose of the terms where $t' \geq t + N^{\eta_2}$.  For this term we observe from \eqref{energy-dispersive} and \ref{ffix-est} that
$$ |Z_N(t',t'')| \lesssim |t'-t''|^{d(\frac{1}{R}-\frac{1}{2})}.$$
Since $d(\frac{1}{R}-\frac{1}{2}) > 2$, the net contribution of these terms is
$$ \lesssim \int_{t' > t + N^{\eta_2}} \int_{t'' < t} |t'-t''|^{d(\frac{1}{R}-\frac{1}{2})}\ dt'' dt'
= O( N^{-\eta_3} ).$$
A similar argument disposes of those terms for which $t'' \leq t - N^{\eta_2}$, so we are left with showing that
$$ |\int_t^{t+N^{\eta_2}} \int_{\max(t-N^{\eta_2},0)}^t Z_N(t',t'')\ dt'' dt'| \lesssim N^{-\eta_3}.$$
By H\"older's inequality we can bound the left-hand side by
\begin{align*}
&\| P_{\geq N} F(u) \|_{L^{q'_0}_t W^{1,r'_0}_x([t,t+N^{\eta_2}]\times\R^d)} \times \\
&\left\| \int_{\max(t-N^{\eta_2},0)}^t e^{i(t'-t'')\Delta} P_{\geq N} F(u(t''))\ dt''
\right\|_{L^{q_0}_t W^{1,r_0}_x([t,t+N^{\eta_2}]\times\R^d)}.
\end{align*}
Applying \eqref{retarded} we may bound the left-hand side by
$$ \lesssim \| P_{\geq N} F(u) \|_{L^{q'_0}_t W^{1,r'_0}_t([\max(t-N^{\eta_2},0), t+N^{\eta_2}])}^2,$$
which by Proposition \ref{smooth-effect} and dyadic decomposition can be bounded by $O( N^{\eta_2/q'_0} N^{- \eta_1} )$, which is acceptable.
\end{proof}

\section{Preliminary spatial localisation}\label{premspac-sec}

In previous sections we decoupled $u$ into a radiation term $e^{it\Delta} u_+$, a weakly bound state $v(t)$, and an asymptotically vanishing error.  The weakly bound state $v$ was localised in frequency, 
but we have not yet achieved (strong) compactness type control on $v$ because we have not localised $v$ in \emph{physical space}.  Of course, in the non-spherically symmetric case the action of the translation group $G$ shows that one cannot hope to naively localise $v(t)$ to be near the origin; the example of multisoliton solutions shows that one cannot even hope to localise $v(t)$ near a single time-varying point $x(t) \in \R^d$.  However, we are still able to localise $v(t)$ to a \emph{bounded number} of time-varying points $x_1(t), \ldots, x_J(t)$.  Indeed, we now assert

\begin{theorem}[Preliminary spatial localisation]\label{premspat}  Let $E > 0$ and $0 < \mu_0 < 1$.  Then there exists $J = J(E,\mu_0)$ and $\mu_5 = \mu_5(E,\mu_0) > 0$ depending only\footnote{Of course, we are implicitly allowing everything to depend also on the fixed parameters $d,p,C_0,\theta$.} on $E, \mu_0$ with the following property: whenever $u$ is a forward-global solution of energy at most $E$, then there exists functions $x_1,\ldots,x_J: \R^+ \to \R^d$ such that
\begin{equation}\label{muloc}
 \limsup_{t \to \infty} \int_{\inf_{1 \leq j \leq J} |x - x_j(t)| \geq 1/\mu_5} |v(t,x)|^2\ dx
\lesssim \mu_0^2,
\end{equation}
where $v$ is the weakly bound component of $u$.
\end{theorem}

\begin{remark} Again, a glance at Proposition \ref{attractive} indicates the relevance of this theorem to Theorems \ref{attract-rad}, \ref{attract-nonrad}.  Note however that the number of concentration points $J$ is currently allowed to depend on the error tolerance $\mu_0$.  This means that this localisation is in fact somewhat weaker than what is necessary to obtain Theorem \ref{attract-nonrad} (though if one assumes spherical symmetry then we will still be able to recover Theorem \ref{attract-rad} without much difficulty).  We shall address this issue of non-uniformity in $J$ in later sections.
\end{remark}

The rest of this section is devoted to the (lengthy) proof of this theorem.  The main tool in this theory shall be the strong decay (better than $t^{-2}$) of the fundamental solution of $e^{it\Delta}$; this decay will be exploited both locally in time (via Strichartz theory) and in the distant future and past (via the Duhamel formula).  

Fix $E > 0$ and $\mu_0$.  In addition to $\mu_0$, we will need some additional
small quantities
$$ \mu_0 \gg \mu_1 \gg \mu_2 \gg \mu_3 \gg \mu_4 > 0$$
depending on $E$ and the other fixed parameters to be chosen later.  Observe that to prove \eqref{muloc} we only need to analyse the behaviour for sufficiently late times $t$, where the definition of ``sufficiently late'' can depend on $E,u$, and the fixed parameters.  We will thus assume that all times are as late as necessary for the arguments which follow.

\subsection{First step: $L^\infty_x$ spatial localisation at fixed times}

Fix a forward-global solution $u$, and let $v$ be the weakly bound component.  The first
step is to exploit the frequency localisation to restrict to medium frequencies.  Define
$v_\med := P_{\mu_2 < \cdot < 1/\mu_2} v$ to be the medium frequency component of $v$.
For sufficiently late times $t$ (depending on all previous parameters), we see from Proposition \ref{asymloc} that
\begin{equation}\label{vvm}
\| v(t) - v_\med(t) \|_{H} \lesssim \mu_2^{\eta_3}
\end{equation}
and thus
\begin{equation}\label{ut}
 u(t) = e^{it\Delta} u_+ + v_\med(t) + O_H(\mu_2^{\eta_3}).
 \end{equation}
Now for each sufficiently late $t$, let $x_1(t), \ldots, x_{J(t)}(t)$ be a maximal $1/\mu_3$-separated set of points in $\R^d$ such that 
\begin{equation}\label{vj}
|v_\med(t,x_j(t))| \geq \mu_3^{1/\eta_1} \hbox{ for all } 1 \leq j \leq J(t).
\end{equation}
From the rapid decay of the convolution kernel of $P_{\mu_2 < \cdot < 1/\mu_2}$ and \eqref{vbound} one can easily establish a bound of the form
$$ |v_\med(t,x_j(t))|^2 \lesssim \int_{|x-x_j(t)| \leq 1/2\mu_3} |v(t,x)|^2\ dx + O( \mu_3^{3/\eta_1} )$$
(say) and so we conclude that
$$ \int_{|x-x_j(t)| \leq 1/2\mu_3} |v(t,x)|^2\ dx \gtrsim \mu_3^{2/\eta_1}$$
for all $j$.  On the other hand, $v(t)$ is uniformly bounded in $H$ norm and thus in $L^2$ norm by $O(1)$.  Thus there exists an integer $J = J(E,\mu_3)$ such that $J(t) \leq J$ for all $t$.  If we then arbitrarily define $x_j(t)$ for $J(t) < j \leq J$ (e.g. setting $x_j(t) = 0$ in these cases) we have thus constructed a sequence $x_1(t),\ldots,x_J(t)$ of points for all sufficiently late times $t$ 
with the property that
\begin{equation}\label{vmedt}
 |v_\med(t,x)| < \mu_3^{1/\eta_1} \hbox{ whenever } \inf_{1 \leq j \leq J} |x-x_j(t)| \geq 1/\mu_3.
\end{equation}

\subsection{Second step: $L^\infty_x$ spatial localisation on a time interval}

We now have $L^\infty_x$ control over a (significant component of) $v(t_0)$ away from the concentration points $x_1(t_0),\ldots,x_J(t_0)$, for any sufficiently late time $t_0$.  Next, we need to use some local theory\footnote{It is possible that the arguments here could be simplified by exploiting mass conservation, as is done in later sections, instead of relying purely on Strichartz theory and approximate finite speed of propagation.  However we feel the arguments here, while lengthier, are more natural, and also have the (minor) benefit of extending to cover non-Hamiltonian nonlinearities.} to also obtain control at times $t$ close to $t_0$.

Fix a sufficiently late time $t_0$, and let $I$ be the interval $I := [t_0-\mu_1^{-1}, t_0+\mu_1^{-1}]$.  Let $D: \R^d \to \R^+$ be the distance function $D(x) := \inf_{1 \leq j \leq J} |x-x_j(t_0)|$; thus \eqref{vmedt} asserts that
$|v_\med(t_0,x)| \leq \mu_3^{1/\eta_1}$ whenever $D \geq 1/\mu_3$.  

Let $\chi: \R^d \to \R^+$ be a smooth cutoff which equals $1$ when $D(x) \leq 2/\mu_3$, vanishes when $D(x) \geq 3/\mu_3$, and obeys the derivative bounds $\nabla^k \chi = O_k( \mu_3^k )$ for $k \geq 0$; such a function can for instance be created by starting
with a smooth function of $D$ which equals $1$ when $D(x) \geq 2.1/\mu_3$ and vanishes when $D(x) \geq 2.9/\mu_3$, and then convolving with an approximation to the identity of width $0.1/\mu_3$; we omit the details.

The portion of $u(t_0)$ away from the concentration points has small linear evolution:

\begin{lemma}\label{late-evolve} We have 
$$
\limsup_{t_0 \to +\infty} \| e^{i(t-t_0)\Delta}[ (1-\chi) u(t_0) ] \|_{L^{q_0}_t W^{1,r_0}_x \cap L^\infty_t L^{Q_0}_x(I \times \R^d)}
\lesssim \mu_2^{\eta_3}.
$$
\end{lemma}

\begin{proof}  
We use the decomposition \eqref{ut}.  We split $u_+$ further as the sum of a test function $\tilde u_+ \in C^\infty_0$ and an error of $H$ norm $O( \mu_2^{\eta_3} )$, leading to the decomposition
$$  (1-\chi) u(t_0) = (1-\chi) e^{it_0\Delta} \tilde u_+ + (1-\chi) v_\med(t_0) + O_H(\mu_2^{\eta_3})$$
(noting the easy fact that multiplication by $\chi$ is bounded on $H$).  We further decompose
\begin{align*}
(1-\chi) v_\med(t_0) &= (1-\chi) P_{<100/\mu_2} (1_{D \geq 1/\mu_3} v_\med(t_0)) \\
&\quad + (1-\chi) P_{<100/\mu_2} (1_{D < 1/\mu_3} v_\med(t_0)).
\end{align*}
From the rapid decay of the convolution kernel of $P_{<100/\mu_2}$, the support of $1-\chi$, and the $H$-boundedness of $v(t_0)$ we see that $(1-\chi) P_{<100/\mu_2} (1_{D < 1/\mu_3} v_\med)$ can be absorbed into the $O_H(\mu_2^{\eta_3})$ error, and we conclude that
\begin{align*}
(1-\chi) u(t_0) &= (1-\chi) e^{it_0\Delta} \tilde u_+ \\
&\quad + (1-\chi) P_{<100/\mu_2} (1_{D \geq 1/\mu_3} v_\med(t_0)) + O_H(\mu_2^{\eta_3}).
\end{align*}
The contribution of $e^{it_0\Delta} \tilde u_+$ is acceptable for sufficiently late times $t_0$ by standard stationary phase estimates (see e.g. \cite{stein:large}), taking advantage of the smoothness and compact support of $\tilde u_+$.  The contribution of the $O_H(\mu_2^{\eta_3})$ error is acceptable by Strichartz estimates \eqref{homog} and Sobolev embedding.  We are thus left to demonstrate that
$$\| e^{i(t-t_0)\Delta}[ (1-\chi) P_{<100/\mu_2} (1_{D \geq 1/\mu_3} v_\med(t_0)) ] \|_{L^{q_0}_t W^{1,r_0}_x \cap L^\infty_t L^{Q_0}_x(I \times \R^d)}
\lesssim \mu_2^{\eta_3}
$$
for sufficiently late times $t_0$.
However, stationary phase (see e.g. \cite{stein:large}) reveals that the operators 
$\nabla^j e^{i(t-t_0)\Delta} (1-\chi) P_{<100/\mu_2}$ have an operator norm of $O_{\mu_1}(\mu_2^{-1/\eta_0})$
(say) on $L^{r_0}_x$ and $L^{Q_0}_x$ for $j=0,1$, so the left-hand side can be bounded by
$$\lesssim_{\mu_1} \mu_2^{-1/\eta_0} 
\| 1_{D \geq 1/\mu_3} v_\med(t_0) \|_{L^{r_0}_x(\R^d) \cap L^{Q_0}_x(\R^d)}.$$
But from \eqref{vmedt} we have
$$\| 1_{D \geq 1/\mu_3} v_\med(t_0) \|_{L^\infty_x(\R^d)} \leq \mu_2^{1/\eta_1}.$$
Interpolating this with \eqref{vbound} we obtain the claim.
\end{proof}

In light of Lemma \ref{late-evolve}, and from the heuristic that medium frequencies propagate at approximately finite speed, we expect the nonlinear solution to also be small away from the points of concentration and for times near $t_0$.  This is indeed the case:

\begin{lemma}[Spatial decay]\label{dmu}  We have
$$ \limsup_{t_0 \to +\infty} \| 1_{D > \mu_4^{-2}} u \|_{L^{q_0}_t L^{r_0}_x(I \times \R^d)} \lesssim_{\mu_1} \mu_2^{\eta_3}.$$
\end{lemma}

\begin{proof} Let $\tilde u$ be the solution to \eqref{nls} on $I$ with initial data $\tilde u(t_0) = \chi u(t_0)$.  From Lemma \ref{late-evolve} and Lemma \ref{perturb} we see that this solution exists on $I$, and in fact we have
\begin{equation}\label{utu}
\| u - \tilde u \|_{C^0_t H^1_x(I \times \R^d)} \lesssim_{\mu_1} \mu_2^{\eta_3}.
\end{equation}
In particular
\begin{equation}\label{tut}
 \tilde u(t) = O_H(1) \hbox{ for all } t \in I.
\end{equation}
We shall need another weight function $W: \R^d \to \R^+$ comparable to $1 + \mu_4 D$ which obeys the derivative bounds $\nabla W, \nabla^2 W = O(\mu_4)$; such a function can for instance be obtained by convolving $1 + \mu_4 D$ with an approximation to the identity at scale $1/\mu_4$.   Since
$W \chi = O(1)$, we have
\begin{equation}\label{wtu}
 \| W \tilde u(t_0) \|_{L^2_x} \lesssim \|u(t_0)\|_{L^2_x} \lesssim 1.
 \end{equation}
On the other hand, since $\tilde u$ obeys \eqref{nls} we see \eqref{fpower} and from the derivative bounds on $W$ that
$$ (i \partial_t + \Delta) (W \tilde u) = O( W |u|^p ) + O( \mu_4 |\tilde u| ) + O( \mu_4 |\nabla \tilde u| ).$$
Applying Strichartz estimates \eqref{homog}, \eqref{retarded} we conclude that\footnote{It is not immediately obvious that the norms here are finite; however we can truncate the weight $W$ to be bounded, apply the estimates here, and then use a monotone convergence argument to remove the truncation.  We omit the details.}
$$ \| W \tilde u \|_{C^0_t L^2_x \cap L^{q_0}_t L^{r_0}_x(I' \times \R^d)} \lesssim \| W \tilde u(t') \|_{L^2_x} + \| W |\tilde u|^p \|_{L^{q'_0}_t L^{r'_0}_x(I' \times \R^d)} + \mu_4 $$
for any sub-interval $I'$ of $I$ and any $t' \in I'$.  If we denote the left-hand side by $X(I')$,
we observe from H\"older, \eqref{tut}, and Sobolev embedding that
$$ \| W |\tilde u|^p \|_{L^{q'_0}_t L^{r'_0}_x(I' \times \R^d)} \lesssim 
|I'|^{1/q'_0-1/q_0} X(I')$$
and hence (for $I'$ sufficiently small) we have that
$$ X(I') \lesssim \| W \tilde u(t') \|_{L^2_x(\R^d)} + \mu_4.$$
Iterating this using \eqref{wtu} (chopping $I$ up into sufficiently small intervals) we conclude that
$$ X(I) \lesssim_{\mu_1} 1;$$
in particular, we see that
$$ \| 1_{D > \mu_4^{-2}} \tilde u \|_{L^{q_0}_t L^{r_0}_x(I \times \R^d)}
 \lesssim_{\mu_1} \mu_4.$$
Combining this with \eqref{utu} we obtain the claim.
\end{proof}

This shows that local-in-time Duhamel effects are localised in space:

\begin{corollary}\label{fnon}  For $t_0$ sufficiently large, and for any $I' \subset I$, we have
$$ \| 1_{D > \mu_4^{-3}} \int_{I'} e^{i(t_0-t')\Delta} F(u(t'))\ dt' \|_{L^2_x(\R^d)} \lesssim_{\mu_1} \mu_2^{\eta_3}.$$
\end{corollary}

\begin{proof} Let $\chi_2: \R^d \to \R^+$ be a smooth cutoff which equals $1$ when $D \leq \mu_4^{-2}$, equals $0$ when $D > 2\mu_4^{-2}$, and obeys the usual derivative estimates in between.
We split
$$F(u(t')) = P_{\geq 1/\mu_2} F(u(t') \chi_2) + P_{\leq 1/\mu_2}[ 1_{D \leq 2\mu_4^{-2}} F(u(t') \chi_2)] + O( 1_{D > \mu_4^{-2}} |u(t')|^p ).$$
From Lemma \ref{dmu}, H\"older, and \eqref{sob} we have
$$ \|1_{D > \mu_4^{-2}} |u(t')|^p\|_{L^{q'_0}_t L^{r'_0}_x(I \times \R^d)} \lesssim_{\mu_1} \mu_2^{\eta_3}$$
so the contribution of the error term $O( 1_{D > \mu_4^{-2}} |u(t')|^p )$ is acceptable from \eqref{dual-homog}.  To deal with the low frequency term $P_{\leq 1/\mu_2}[ 1_{D \leq 2\mu_4^{-2}} F(u(t') \chi)]$, we observe from Schur's test\footnote{Schur's test only establishes $L^2 \to L^2$ boundedness, but one can establish $L^1 \to L^2$ bounds by Minkowski's inequality and then interpolate.} and stationary phase (see e.g. \cite{stein:large}) that the operator
$1_{D > \mu_4^{-3}} e^{i(t_0-t')\Delta} P_{\leq 1/\mu_2} 1_{D \leq 2\mu_4^{-2}}$ has
an $L^R_x \to L^2_x$ operator norm of $O( \mu_4^{100} )$ (say) for each $t \in I$, 
while from \eqref{ffix-est} we have
$\| F(u(t')) \chi_2 \|_{L^R_x} = O(1)$.  So this term is also acceptable.

Finally, the contribution of the high frequency term $P_{\geq 1/\mu_2} F(u(t') \chi_2)$ can be
controlled using \eqref{dual-homog} by
$$ \lesssim \| P_{\geq 1/\mu_2} F(u \chi) \|_{L^{q'_0}_t L^{r'_0}_x(I \times \R^d)}
\lesssim \mu_2 \| \nabla F(u \chi_2) \|_{L^{q'_0}_t L^{r'_0}_x(I \times \R^d)}.$$
But this is acceptable by a minor variant of \eqref{naff} (one can easily repeat the proof of this estimate and see that the smooth cutoff $\chi_2$ causes no additional difficulty).
\end{proof}

\subsection{Third step: $L^2_x$ localisation at fixed times}

We can now upgrade our $L^\infty$ localisation of $v$ to $L^2$ localisation:

\begin{proposition}\label{vl2}
Let $\chi_3: \R^d \to \R^+$ be a smooth cutoff which equals $1$ when $D \geq \mu_4^{-3}$, equals $0$ when $D \leq 2\mu_4^{-3}$, and obeys the usual derivative estimates in between.  Then for $t_0$ sufficiently large, we have
$$ \chi_3 v(t_0) = O_{L^2}(\mu_1^c)$$
for some $c > 0$ (depending only on $d,p$).
\end{proposition}

\begin{proof}
From Corollary \ref{fnon} we see that
$$ \chi_3 \int_{I'} e^{i(t_0-t')\Delta} F(u(t'))\ dt' = O_{L^2}(\mu_2^{\eta_3/2})$$
(say) for all $I' \subset I$.  In particular from Duhamel's formula we have
$$ \chi_3 v(t_0) = \chi_3 e^{- i\mu_1^{-1} \Delta} v(t_0 + 1/\mu_1) + O_{L^2}(\mu_2^{\eta_3/2})$$
and
$$ \chi_3 v(t_0) = \chi_3 e^{+ i\mu_1^{-1} \Delta} v(t_0 - 1/\mu_1) + O_{L^2}(\mu_2^{\eta_3/2}).$$
Taking the inner product of these two estimates (and using \eqref{energy}) we conclude that
$$ \| \chi_3 v(t_0)\|_{L^2_x(\R^d)}^2
= \langle \chi_3 e^{- i\mu_1^{-1} \Delta} v(t_0 + 1/\mu_1), \chi_3 e^{i\mu_1^{-1} \Delta} v(t_0 - 1/\mu_1) \rangle_{L^2_x} + O( \mu_2^{\eta_3/2} ).$$
Combining the two $\chi_3$ factors together, it thus suffices to show that
$$ |\langle e^{- i\mu_1^{-1} \Delta} v(t_0 + 1/\mu_1), \chi_3^2 e^{i\mu_1^{-1} \Delta} v(t_0 - 1/\mu_1) \rangle_{L^2_x}| \lesssim \mu_1^{2c}.$$
We approximate $u(0)-u_+$ as $u(0)-u_+ = \phi + O_{L^2}(\mu_2)$ (say), where $\phi$ is a Schwartz function (independent of $t_0$), thus by
\eqref{uforward}
$$ v(t_0-1/\mu_1) = e^{i(t_0-1/\mu_1)\Delta} \phi - i \int_0^{t_0-1/\mu_1} e^{i(t_0-1/\mu_1-t'')\Delta} F(u(t''))\ dt'' + O_{L^2}(\mu_2)$$
and hence
\begin{align*} \chi_3^2 e^{i\mu_1^{-1} \Delta} v(t_0 - 1/\mu_1) &=
e^{it_0\Delta} \phi - (1-\chi_3^2) e^{it_0\Delta} \phi \\
&- i \chi_3^2 \int_0^{t_0-1/\mu_1} e^{i(t_0-t'')\Delta} F(u(t''))\ dt'' + O_{L^2}(\mu_2).
\end{align*}
From dispersive estimates we see that $(1-\chi_3^2) e^{it_0\Delta} \phi = O_{L^2}(\mu_2)$ for sufficiently large $t_0$,
so that term may be absorbed into the error term.  The contribution of the $O_{L^2}(\mu_2)$ error can be discarded by Cauchy-Schwarz.  Using this and \eqref{ubackward}, we thus reduce to showing that
\begin{equation}\label{mu1-1} \int_{t_0+1/\mu_1}^{+\infty} 
|\langle e^{i(t_0-t')\Delta} F(u(t')), e^{it_0\Delta} \phi \rangle_{L^2_x}|\ dt' \lesssim \mu_1^{2c}
\end{equation}
and
\begin{equation}\label{mu1-2}
 \int_{t_0+1/\mu_1}^{+\infty} \int_0^{t_0-1/\mu_1}
|\langle e^{i(t_0-t')\Delta} F(u(t')), \chi^2_3 e^{i(t_0-t'')\Delta} F(u(t'')) \rangle_{L^2_x}|\ dt' dt'' \lesssim \mu_1^{2c}.
\end{equation}
Let us first prove \eqref{mu1-1}.  By \eqref{ffix-est} and H\"older we have
$$ |\langle e^{i(t_0-t')\Delta} F(u(t')), e^{it_0\Delta} \phi \rangle_{L^2_x}|
= |\langle F(u(t')), e^{it'\Delta} \phi \rangle_{L^2_x}| \lesssim \| e^{it'\Delta} \phi\|_{L^{R'}_x(\R^d)}$$
but by dispersive estimates we have $\| e^{it'\Delta} \phi\|_{L^{R'}_x(\R^d)} \ll_\phi (t')^{-d(\frac{1}{R}-\frac{1}{2})}$.  By choice of $R$ we have $d(\frac{1}{R}-\frac{1}{2}) > 2$, and so we obtain \eqref{mu1-1} for $t_0$ sufficiently large.

Now we prove \eqref{mu1-2}.  Writing $\chi_3^2 = 1 - (1-\chi_3^2)$ and using \eqref{prop} we have
$$ \langle e^{i(t_0-t')\Delta} f, \chi^2_3 e^{i(t_0-t'')\Delta} g \rangle_{L^2_x}
= \langle e^{i(t''-t')\Delta} f, g \rangle_{L^2_x} -
\int_{\R^d} \int_{\R^d} f(x) \overline{g(z)} K_{t',t''}(x,z)\ dx dz$$
for arbitrary test functions $g$, where $K_{t',t''}$ is the kernel
\begin{align*}
 K_{t',t''}(x,z) &:= \frac{1}{(4\pi i (t_0-t'))^{d/2} (-4\pi i (t_0-t''))^{d/2}} \\
&\quad \int_{\R^d} e^{i |x-y|^2 / 4(t_0-t')} e^{-i|x-y|^2/4(t_0-t'')} (1-\chi^2_3)(y)\ dy.
\end{align*}
The $|y|^2$ coefficient of the phase in $K_{t',t''}$ is $O( \frac{|t'-t''|}{|t_0-t'| |t_0-t''|} )$.  Applying stationary phase (see e.g. \cite{stein:large}) one concludes that
$$ |K_{t',t''}(x,z)| \lesssim |t'-t''|^{-d/2}.$$
(Note that each derivative of $1-\chi^2_3$ picks up some powers of $\mu_4$, which eventually compensate for the measure of the support of this function, which is $O( J (\mu_4^3)^{-d})$, since $J$ is bounded by a quantity depending only on $\mu_3$.  Thus the contribution of the regions of space where the phase is non-stationary can be dealt with by repeated integration by parts in the usual fashion.)
From this and \eqref{prop} we conclude that
$$ \langle e^{i(t_0-t')\Delta} f, \chi^2_3 e^{i(t_0-t'')\Delta} g \rangle = O( |t'-t''|^{-d/2} \|f\|_{L^1_x(\R^d)}
\|g\|_{L^1_x(\R^d)} ).$$
On the other hand, from Cauchy-Schwarz we have
$$ \langle e^{i(t_0-t')\Delta} f, \chi^2_3 e^{i(t_0-t'')\Delta} g \rangle = O( \|f\|_{L^2_x(\R^d)}
\|g\|_{L^2_x(\R^d)} )$$
and hence by bilinear interpolation
$$ \langle e^{i(t_0-t')\Delta} f, \chi^2_3 e^{i(t_0-t'')\Delta} g \rangle = O( |t'-t''|^{-d(\frac{1}{2}-\frac{1}{R})} \|f\|_{L^R_x(\R^d)} \|g\|_{L^R_x(\R^d)} )$$
and thus by \eqref{ffix-est}
$$ \langle e^{i(t_0-t')\Delta} F(u(t')), \chi^2_3 e^{i(t_0-t'')\Delta} F(u(t'')) \rangle
= O( |t'-t''|^{-d(\frac{1}{2}-\frac{1}{R})} ).$$
Again, by choice of $R$ we have $d(\frac{1}{2}-\frac{1}{R}) > 2$, and \eqref{mu1-2} follows (for $c$ sufficiently small).
\end{proof}

The claim \eqref{muloc} is now immediate from Proposition \ref{vl2}.  This proves Theorem \ref{premspat} (defining $\mu_5$ appropriately).

\begin{remark} It is possible (basically thanks to \eqref{pnbig}) to obtain a similar result to \eqref{muloc} with the mass density $|v(t,x)|^2$ replaced by the energy density $|v(t,x)|^2 + |\nabla v(t,x)|^2$.  We will not need this apparently stronger localisation in our argument, though, and in any event it will eventually follow from our precompactness theorems thanks to Proposition \ref{precompac}.
\end{remark}

\section{The radial case}\label{radial-sec}

We now have enough control on forward-global solutions to obtain the desired compactness in the spherically symmetric case.  Strictly speaking, the arguments in this section are redundant, as they will be superceded by the more general non-radial arguments in later sections, but we present them here to offer a simplified version of the arguments to come.

We first note that the spherical symmetry allows one to collapse the $J$ points of concentration to a single one, namely the origin:

\begin{proposition}[Radial spatial localisation]\label{radspac}  Let $E > 0$ and $0 < \mu_0 < 1$.  Then there exists $\mu_6 = \mu_6(E,\mu_0)$ depending only on $E, \mu_0$ (and on $d,p,\theta$) with the following property: whenever $u$ is a spherically symmetric forward-global solution of energy at most $E$, we have
\begin{equation}\label{muloc-spin}
\limsup_{t \to +\infty} \int_{|x| \geq 1/\mu_6} |v(t,x)|^2 \ dx \lesssim \mu_0^2,
\end{equation}
where $v$ is the weakly bound component of $u$.
\end{proposition}

\begin{proof}  Let $E, \mu_0, u$ be as above.  We apply Theorem \ref{premspat} to conclude that for all sufficiently large times $t$ we have $x_1(t),\ldots,x_J(t)$ such that
$$ \int_{\R^d} 1_{\inf_{1 \leq j \leq J} |x - x_j(t)| \geq 1/\mu_5} |v(t,x)|^2\ dx
\lesssim \mu_0^2.$$
Since $u$ is spherically symmetric, it is easy to see (from uniqueness and rotational symmetry) that $u_+$ and $v$ are also spherically symmetric.  Thus we may average the above estimate over rotations and conclude that
$$ \int_{\R^d} 
(\int_{S^{d-1}} 1_{\inf_{1 \leq j \leq J} ||x|\omega - x_j(t)| \geq 1/\mu_5}\ d\omega)
|v(t,x)|^2\ dx
\lesssim \mu_0^2.$$
In particular, we have
$$ \int_{|x| \geq 1/\mu_6} 
(\int_{S^{d-1}} (1 - \sum_{j=1}^J 1_{||x|\omega - x_j(t)| < 1/\mu_5})\ d\omega)
|v(t,x)|^2\ dx
\lesssim \mu_0^2.$$
But if $\mu_6$ is sufficiently small depending on $\mu_5$, we see from elementary trigonometry that
$$ \int_{S^{d-1}} 1_{||x|\omega - x_j(t)| < 1/\mu_5}\ d\omega \lesssim (\mu_6/\mu_5)^{d-1}$$
so if $\mu_6$ is small enough depending on both $J$ and $\mu_5$, we have
$$ \int_{S^{d-1}} (1 - \sum_{j=1}^J 1_{||x|\omega - x_j(t)| < 1/\mu_5})\ d\omega \geq \frac{1}{2}$$
for all $|x| \geq 1/\mu_6$, and the claim follows.
\end{proof}

We can now quickly prove Theorem \ref{attract-rad} and Corollary \ref{petite-rad}.

\begin{proof}[Proof of Theorem \ref{attract-rad}]  From Proposition \ref{asymloc}, Proposition \ref{radspac}, and Proposition \ref{attractive} we see that there exists a compact set $K \subset H$ such that $\lim_{t \to +\infty} \dist_H(u(t),K)=0$ for all spherically symmetric forward-global solutions $u$ of energy at most $E$.

At present, $K$ is not necessarily invariant under the NLS flow.  To address this, let $\K_{E,\rad}$ be the closure of the set of all limit points $\lim_{n \to \infty} v(t_n)$ of weakly bound states associated to forward-global spherically symmetric solutions of energy $E$, where $t_n$ ranges over sequences of times that go to infinity.  This is clearly a closed and hence compact subset of $K$, and is also clearly an attractor for $v(t)$ (again thanks to the compactness of $K$).  The local well-posedness theory also ensures that this set is invariant under $S(t)$ for all small $t$, and hence for all large $t$ also, and we are done.

Finally we address the uniqueness of $u_+$.  If we had two $u_+, \tilde u_+$ obeying \eqref{tsim}, then on subtraction we see that the set $\{ e^{it\Delta} (u_+ - \tilde u_+): t \geq 0 \}$ is totally bounded, and hence precompact.  Combining this with Lemma \ref{rl} we see that $e^{it\Delta}(u_+ - \tilde u_+)$ converges to zero weakly (and hence in $H$ norm, by precompactness) as $t \to \infty$.  But $e^{it\Delta}$ is unitary, and hence $u_+ = \tilde u_+$, a contradiction.
\end{proof}

\begin{proof}[Proof of Corollary \ref{petite-rad}]  From Proposition \ref{precompac} we see that (i) implies (iv), which trivially implies (iii).  
Next we show that (iii) implies (ii).  It suffices to show that $\langle u_+, \phi \rangle = 0$ for every test function $\phi \in C^\infty_0(\R^d)$.  From \eqref{uplus} and unitarity it suffices to show that
$$ \limsup_{t \to +\infty} |\langle u(t), e^{it\Delta} \phi \rangle| = 0.$$
But from Lemma \ref{rl} we have
$$ \limsup_{t \to \infty} \int_{|x| \leq R} |u(t,x)| |e^{it\Delta} \phi(x)|\ dx = 0$$
for all $R > 0$, while from the hypothesis (iii) and Cauchy-Schwarz we have
$$ \limsup_{t \to \infty} \int_{|x| > R} |u(t,x)| |e^{it\Delta} \phi(x)|\ dx \to 0 \hbox{ as } R \to \infty.$$
The claim follows.

Finally we show that (ii) implies (i).  Since $u_+ = 0$, we have $u = v$, and hence by Theorem \ref{attract-rad}
we have $\dist_H(u(t),\K_{E,\rad}) \to 0$ as $t \to +\infty$.  Since $\K_{E,\rad}$ was compact, this (combined with the continuity of the NLS flow) implies that the orbit $\{ u(t): t \geq 0 \}$ is totally bounded, and thus $u$ is almost periodic as claimed.
\end{proof}

\section{Final spatial localisation}\label{fsl-sec}

We now return to the non-radial case.  From Propositions \ref{pre}, \ref{asymloc} and Theorem \ref{premspat} we have
already achieved much of Theorem \ref{attract-nonrad}.  The main remaining technical issue is that the number of concentration points $J(E,\mu_0)$ in Theorem \ref{premspat} can go to infinity as $\mu_0$ goes to zero, which prohibits us from obtaining an attractor which is precompact in the sense of Definition \ref{gpre}.  To prevent this, we need to show that each concentration point in fact absorbs a large portion of the mass of $u$.  More precisely, we shall establish the following estimate:

\begin{theorem}[Mass concentration property]\label{radspac-2}  Let $u$ be a forward-global solution of energy at most $E$.  Suppose that we have the mass concentration bound
$$ \int_{|x-x_0| < R} |u(t_0,x)|^2\ dx \geq \mu_1^2$$
for some $x_0 \in \R^d$, $t_0 \in \R^+$, $R > 0$, and $\mu_1 > 0$.  Then, if $t_0$ is sufficiently large depending on $u,E,x_0,R,\mu_1$, we have the improved mass concentration bound
$$ \int_{|x-x_0| < R'} |u(t_0,x)|^2\ dx \gg 1$$
for some $R' = R'(E,R,\mu_1) < \infty$ depending only on $E,R,\mu_1$ (in particular, $R'$ is independent of $u$).
\end{theorem}

The rest of this section is devoted to the proof of this theorem.  We first observe that it suffices to prove a weaker statement in which only a slight mass improvement is obtained:

\begin{proposition}[Mass concentration property, inductive step]\label{radspac-3}  Given any $E > 0$ there exists $\mu_0 = \mu_0(E) > 0$ with the following property: Suppose that we have the mass concentration bound
$$ \mu_1^2 \leq \int_{|x-x_0| < R} |u(t_0,x)|^2\ dx \leq \mu_0^2$$
for some $x_0 \in \R^d$, $t_0 \in \R^+$, $R > 0$, and $\mu_1 > 0$, and some forward-global solution $u$ of energy at most $E$.  Then, if $t_0$ is sufficiently large depending on $u,E,x_0,R,\mu_1$, we have the improved mass concentration bound
$$ \int_{|x-x_0| < R'} |u(t_0,x)|^2\ dx \geq \int_{|x-x_0| < R} |u(t_0,x)|^2\ dx + \mu_4^2$$
for some $\mu_4 = \mu_4(E,\mu_1) > 0$ depending only on $E,\mu_1$ and $R' = R'(E,R,\mu_1,\mu_4) < \infty$ depending only on $E,R,\mu_1,\mu_4$.
\end{proposition}

Indeed, Theorem \ref{radspac-2} easily follows from iterating Proposition \ref{radspac-3} at most $\mu_0^2/\mu_4^2$ times.

We now prove Proposition \ref{radspac-3}.  By translation invariance we can take $x_0=0$.
Let $E > 0$, and let $\mu_0 > 0$ be a sufficiently small quantity to be chosen later.  Let $\mu_1, x_0 = 0, t_0, R, u$ be as in the Proposition.  We need some further small quantities
$$ \mu_1 \gg \mu_2 \gg \mu_3 \gg \mu_4 > 0$$
with each $\mu_i$ assumed sufficiently small depending on the previous $\mu_j$ and on $E$.  Finally we let $R' = R'(E,\mu_0,\ldots,\mu_4) > 1$ be a large radius to be chosen later.  In particular we may take $R' > 100 R$.    Let $I$ be the time interval $I := [t_0,t_0+1/\mu_3]$.

Suppose for contradiction that the claim failed, then we have very little mass in an annulus at time $t_0$:
$$ \int_{R < |x| < R'} |u(t_0,x)|^2\ dx \leq \mu_4^2.$$
In particular
$$ \int_{|x| < R'} |u(t_0,x)|^2\ dx \lesssim \mu_0^2.$$
The idea now is to localise the small data scattering theory to the cylinder $\{ (t,x): t \in I; |x| \leq R' \}$, and use this to contradict the spatial localisation from Theorem \ref{premspat}.  We first
use local conservation of mass to extend the above estimates to the rest of the time interval $I$:

\begin{lemma}[Local absence of mass]  If $R'$ is sufficiently large depending on $E, \mu_0, \ldots,\mu_4$, we have
\begin{equation}\label{mu2}
 \sup_{t \in I} \int_{R'/16 < |x| < R'/2} |u(t,x)|^2\ dx \lesssim \mu_4^2
 \end{equation}
and
\begin{equation}\label{muo}
\sup_{t \in I} \int_{|x| < R'/2} |u(t,x)|^2\ dx \lesssim \mu_0^2
\end{equation}
and
\begin{equation}\label{mu-low}
\inf_{t \in I} \int_{|x| < R'/16} |u(t,x)|^2\ dx \gtrsim \mu_1^2.
\end{equation}
\end{lemma}

\begin{proof}  
Let $\varphi:\R^d \to \R^+$ be a smooth non-negative cutoff function supported on the annulus $\{ R'/32 \leq |x| \leq R' \}$ which
equals one on the annulus $\{ R'/16 \leq |x| \leq R'/2 \}$.  Then we have
$$ \int_{\R^d} \varphi^2(x) |u(t_0,x)|^2\ dx \lesssim \mu_4^2.$$
From the Hamiltonian nature of NLS \eqref{nls} we have the formal identity
$$ \partial_t |u|^2 = - 2 \nabla \Im( \overline{u} \nabla u )$$
and thus by Stokes theorem
$$ \int_{\R^d} \varphi^2(x) |u(t,x)|^2\ dx - \int_{\R^d} \varphi^2(x) |u(t_0,x)|^2\ dx
= 4 \int_{t_0}^{t}
\int_{\R^d} \varphi(x) \nabla \varphi(x) \cdot \Im( \overline{u} \nabla u )(t',x)\ dx dt'$$
for all $t \in I$.
This identity can be justified first for smooth $u$ and smooth Hamiltonian nonlinearities $F$, and then by limiting arguments to energy-class $u$ and then to general Hamiltonian nonlinearities $F$.  Applying
Cauchy-Schwarz, \eqref{energy} and the bound $\nabla \varphi = O(1/R')$ we see that
$$ \int_{\R^d} \varphi^2(x) |u(t,x)|^2\ dx \lesssim 
\mu_4 + \frac{1/\mu_3}{R'} \sup_{t' \in I} \left(\int_{\R^d} \varphi^2(x) |u(t',x)|^2\ dx\right)^{1/2}.$$
Taking suprema in $t'$ and assuming $R'$ large depending on $\mu_3,\mu_4$ we obtain
$$ \sup_{t \in I} \int_{\R^d} \varphi^2(x) |u(t,x)|^2\ dx \lesssim 
\mu_4 $$
and the first claim follows.  The other two claims \eqref{muo}, \eqref{mu-low} are proven similarly.
\end{proof}

Now take $\chi:\R^d \to \R^+$ be a smooth cutoff supported on the ball $\{ |x| \leq R'/4\}$ which equals
one on $\{|x| \leq R'/8\}$, and let $w(t) := u(t) \chi$.  Then from \eqref{nls} we see that $w$ solves the forced NLS equation
\begin{align*}
iw_t + \Delta w &= F(w) + 2 \nabla \chi \cdot \nabla u + (\Delta \chi) u + [ F(u) \chi - F(u\chi)]\\
&= O(|w|^p) + O( \frac{1}{R'} (|\nabla u| + |u|) + 1_{R'/8 \leq |x| \leq R'/4} O( |u|^p ).
\end{align*}
By Strichartz \eqref{homog}, \eqref{retarded} we conclude that
\begin{align*}
\|w\|_{L^2_t L^{2d/(d-2)}_x(I \times \R^d)}
&\lesssim \| w(t_0) \|_{L^2_x(\R^d)} + \| |w|^p \|_{L^2_t L^{2d/(d+2)}_x(I \times \R^d)}\\
&\quad + \frac{1}{R'} \| u \|_{L^1_t H^1_x(I \times \R^d)}
+ \| 1_{R'/8 \leq |x| \leq R'/4} |u|^p \|_{L^2_t L^{2d/(d+2)}_x(I \times \R^d)}.
\end{align*}
For brevity we now suppress the $I \times \R^d$ domain.  From \eqref{muo} we see that
$\|w(t_0)\|_{L^2_x(\R^d)} = O(\mu_0)$.
By H\"older's inequality and the hypotheses $1 + \frac{4}{d} < p < 1 + \frac{4}{d-2}$ we have
$$ \| |w|^p \|_{L^2_t L^{2d/(d+2)}_x} 
\leq \| w \|_{L^2_t L^{2d/(d-2)}_x} \|w\|_{L^\infty_t L^2_x}^\sigma \|w\|_{L^\infty_t L^{2d/(d-2)}_x}^{1-\sigma}$$
for some $0 < \sigma < 1$ depending only on $d,p$; from \eqref{energy}, \eqref{muo} we conclude that
$$ \| |w|^p \|_{L^2_t L^{2d/(d+2)}_x}  \lesssim \mu_0^\sigma \| w \|_{L^2_t L^{2d/(d-2)}_x}.$$
Meanwhile, from \eqref{energy} we have
$$ \| u \|_{L^1_t(I \to H)} \lesssim |I| = 1/\eta_3 $$
and from H\"older we have
$$ \| 1_{R'/8 \leq |x| \leq R'/4} |u|^p \|_{L^2_t L^{2d/(d+2)}_x}
\leq \| u \|_{L^2_t L^{2d/(d-2)}_x} \| 1_{R'/8 \leq |x| \leq R'/4} u \|_{L^\infty_t L^2_x}^\sigma \|u\|_{L^\infty_t L^{2d/(d-2)}_x}^{1-\sigma}$$
and hence by Lemma \ref{localst}, \eqref{energy}, \eqref{mu2}
$$ \| 1_{R'/8 \leq |x| \leq R'/4} |u|^p \|_{L^2_t L^{2d/(d+2)}_x} \lesssim |I|^{1/2} \mu_4^\sigma
= \mu_3^{-1/2} \mu_4^\sigma.$$
Putting this all together we obtain
$$ \|w\|_{L^2_t L^{2d/(d-2)}_x}
\lesssim \mu_0 + \mu_0^\sigma \|w\|_{L^2_t L^{2d/(d-2)}_x} + (R')^{-1} +
\mu_3^{-1/2} \mu_4^\sigma.$$
Thus if the $\mu_i$ are chosen sufficiently small (depending on all previous $\mu_j$) and $R'$ is sufficiently large (depending on all the $\mu_i$) we conclude that
$$ \|w\|_{L^2_t L^{2d/(d-2)}_x} \lesssim \mu_0.$$
In particular
$$
\|1_{|x| \leq R'/16} u \|_{L^2_t L^{2d/(d-2)}_x} \lesssim \mu_0.
$$
Applying Lemma \ref{rl}, we conclude (if $t_0$ is sufficiently large depending on $\mu_0, T, u$) the spacetime smallness bound
$$
\|1_{|x| \leq R'/16} v \|_{L^2_t L^{2d/(d-2)}_x} \lesssim \mu_0.
$$
This is a localised analogue of the usual Strichartz-type scattering bounds for small data solutions.
By the pigeonhole principle, we can thus find $t \in I$ such that
\begin{equation}\label{mu3}
\|1_{|x| \leq R'/16} v(t) \|_{L^{2d/(d-2)}_x(\R^d)} \lesssim \mu_0 / |I|^2 \lesssim \mu_3^{1/2}.
\end{equation}

Fix this $t$.  From \eqref{mu-low} and Lemma \ref{rl} we see (taking $t$ sufficiently large depending on $\mu_1, u$) that
$$ \int_{|x| < R'/16} |v(t,x)|^2\ dx \gtrsim \mu_1^2  $$
Also, from Proposition \ref{premspat}, if $t_0$ is sufficiently large depending on $\mu_1, E$, and $\mu_2$ is sufficiently small depending on $\mu_1, E$, then we can find $J = J(E,\mu_1)$ and $x_1(t),\ldots,x_J(t)$ such that
$$
 \int_{\inf_{1 \leq j \leq J} |x - x_j(t)| \geq 1/\mu_2} |v(t,x)|^2\ dx
\lesssim \mu_1^3
$$
(say). Subtracting, we conclude
$$ \int_{|x| < R'/16} 1_{\inf_{1 \leq j \leq J} |x - x_j(t)| < 1/\mu_2}
|v(t,x)|^2\ dx \gtrsim \mu_1^2. $$
The set $\{x: \inf_{1 \leq j \leq J} |x - x_j(t)| < 1/\mu_2\}$ has volume at most $O( J / \mu_2^d ) = O_{\mu_1,\mu_2}(1)$.  By H\"older we conclude that
$$ 
\| 1_{|x| \leq R'/16} v(t) \|_{L^{2d/(d-2)}_x(I \times \R^d)} \gtrsim_{\mu_1,\mu_2} 1$$
which contradicts \eqref{mu3}.  This proves Proposition \ref{radspac-3}, and Theorem \ref{radspac-2} follows.

\section{Constructing the $G$-compact attractor}\label{nonradial-sec}

In this section we use the above localisations to 
establish Theorem \ref{attract-nonrad} and Corollaries \ref{pet1}, \ref{pet2}.
We begin with a non-radial counterpart of Proposition \ref{radspac}:

\begin{proposition}[Non-radial spatial localisation]\label{nonradspac}  
Let $u$ be a forward-global solution of energy at most $E > 0$.  Then there exists $J = J(E) > 0$ depending only on $E$, and functions $x_1,\ldots,x_J: \R^+ \to \R^d$, such that we have following the asymptotic spatial localisation property: given any $\mu_2 > 0$ there exists $\mu_4 = \mu_4(E,\mu_2) > 0$ depending on $E$ and $\mu_2$ (but \emph{independent} of $u$) such that
\begin{equation}\label{muloc-spin-2}
\limsup_{t \to +\infty} 
\int_{\R^d} 1_{\inf_{1 \leq j \leq J} |x - x_j(t)| > 1/\mu_4} |v(t,x)|^2\ dx \lesssim \mu_2^2
\end{equation}
where $v$ is the weakly bound component of $u$.
\end{proposition}

Note that this proposition implies a strengthened version of Theorem \ref{premspat} in which $J$ is  independent of $\mu_0$.

\begin{proof}  Fix $E$.  Let $1 \gg \mu_0 \gg \mu_1 > 0$ be small parameters, with $\mu_0$ sufficiently small depending on $E$ and $\mu_1$ sufficiently small depending on $\mu_0,E$.  Let $u$ be a forward-global solution of energy at most $E$.  Applying Theorem \ref{premspat} (and assuming $\mu_1$ small enough) we can find $J = J(E,\mu_0) > 0$ and functions $x_1,\ldots,x_J: \R^+ \to \R^d$ such that
\begin{equation}\label{jimbo}
\int_{D \geq 1/\mu_1} |v(t,x)|^2\ dx
\lesssim \mu_0^2
\end{equation}
for all sufficiently late times $t$, where $D = D(t,x)$ is the quantity
$$ D := \inf_{1 \leq j \leq J} |x - x_j(t)|.$$  

To prove \eqref{muloc-spin-2} we may restrict attention to small $\mu_2$; in particular we may assume $\mu_2$ sufficiently small depending on $E,\mu_0,\mu_1$.  Let $\mu_3 > 0$ be sufficiently small depending on $E,\mu_2$, and let $\mu_4 > 0$ be sufficiently small depending on $E,\mu_3$.  Applying Theorem \ref{premspat} once again (and assuming $\mu_3$ small enough) we can find $J' = J'(E,\mu_2) > 0$
and $y_1,\ldots,y_{J'}: \R^+ \to \R^d$ such that
\begin{equation}\label{info}
 \int_{\inf_{1 \leq j' \leq J'} |x - y_{j'}(t)| \geq 1/\mu_3} |v(t,x)|^2\ dx
\lesssim \mu_2^2
\end{equation}
for all sufficiently late times $t$.

The next step is to relate the $x_j(t)$ to the $y_{j'}(t)$.  To do this, let us first assume that $t$ is a sufficiently late time and $1 \leq j' \leq J'$ is such that
\begin{equation}\label{lo}
\int_{|x - y_{j'}(t)| < 1/\mu_3} |v(t,x)|^2 \gtrsim \mu_2^2 / J'.
\end{equation}
Applying Theorem \ref{radspac-2} (and assuming $\mu_4$ small enough), we conclude that
$$ \int_{|x - y_{j'}(t)| < 1/2\mu_4} |v(t,x)|^2 \gtrsim 1.$$
Comparing this with \eqref{jimbo} we see (if $\mu_0$ is chosen sufficiently small) that 
$$ \{ x: |x - y_{j'}(t)| < 1/2\mu_4 \} \cap \{ x: \inf_{1 \leq j \leq J} |x - x_j(t)| \geq 1/\mu_1 \}
\neq \emptyset$$
and thus by the triangle inequality we have $|y_{j'}(t) - x_j(t)| \leq 1/2\mu_4 + 1/\mu_1$ for some $j$.  In particular we have
$$ \{ x: |x - y_{j'}(t)| < 1/\mu_3 \} \subset \{ x: D < 1/\mu_4 \}$$
whenever \eqref{lo} holds. In particular, regardless of whether \eqref{lo} holds or not, we have
$$ \int_{|x - y_{j'}(t)| < 1/\mu_3} 1_{D \geq 1/\mu_4}
|v(t,x)|^2 \lesssim \mu_2^2 / J'$$
for all $1 \leq j' \leq J'$.  Summing in $j'$ and using \eqref{info} we conclude that
$$
\int_{D \geq 1/\mu_4} |v(t,x)|^2\ dx
\lesssim \mu_2^2$$
and the claim follows.
\end{proof}

\begin{proof}[Proof of Theorem \ref{attract-nonrad}]
From Proposition \ref{asymloc}, Proposition \ref{nonradspac}, and Proposition \ref{G-attractive} we can locate a compact set $K$ such that
\begin{equation}\label{jgk}
\lim_{t \to +\infty} \dist_H(v(t), J(GK)) = 0
\end{equation}
for all forward-global solutions $u$ of energy at most $E$.  By adding $0$ to $K$ if necessary we may assume that $0 \in K$, and thus $J_1(GK) \subset J_2(GK)$ whenever $J_1 < J_2$.

Next, let $\K_E$ be the subset of $J(GK)$ consisting of those $f \in J(GK)$ such that $S(t) f$ is well-defined and lies in $J(GK)$ for all $t \in \R$.  This is clearly a flow-invariant subset of $J(GK)$ and thus $G$-precompact with $J$ components; it is also $G$-invariant (i.e. translation-invariant) thanks to the $G$-invariance of the NLS equation and of the set $J(GK)$.  From Corollary \ref{compact-k} and the continuous nature of the flow maps $S(t)$ we see that $\K_E$ is closed.

Now we establish the profile decomposition \eqref{profdecomp} (which implies \eqref{tsim-nonrad} as can be seen by a proof by contradiction).  Let $u$ be a forward-global solution of energy at most $E$, and let $t_n$ be a sequence of times going to infinity.  From \eqref{jgk} and Lemma \ref{cc} we can (after passing to a subsequence) obtain a representation
$$ v(t_n) = \sum_{m=1}^M \tau_{x_{m,n}} w_m + o_H(1)$$
for some $J_1 + \ldots + J_M = J$, $w_m \in J_m(GK)$, and $x_{m,n}$ obeying the asymptotic separation condition \eqref{ass}.

Now let $\mu_0 > 0$ be a sufficiently small quantity (depending only on $E$).  Suppose that $t$ is a time with $|t| \leq \mu_0$.  By Lemma \ref{approx} and Lemma \ref{asl}, we have
$$ v(t_n+t) = \sum_{m=1}^M \tau_{x_{m,n}} S(t) w_m + o_H(1)$$
for sufficiently large $n$.  Applying Lemma \ref{cc-inverse} we conclude that there exists a partition
$J = J_1(t) + \ldots + J_M(t)$ such that $S(t) w_m \in J_m(t)(GK)$.  In particular $S(t) w_m \in J(GK)$ for all $t$, which implies that $w_m \in \K_E$.  This gives \eqref{profdecomp} as desired.

Finally, we show that the radiation state $u_+$ is unique.  If this were not the case, then we could find distinct $u_+, \tilde u_+$ obeying \eqref{tsim-nonrad}.  But then we have
$$ u_+ - \tilde u_+ = e^{-it\Delta} (w(t) - \tilde w(t)) + o_H(1)$$
for some $w(t), \tilde w(t) \in J\K_E$.  Using Corollary \ref{grl} we conclude that $\|u_+ - \tilde u_+ \|_{L^q(\R^d)} = 0$ for all $2 < q \leq \frac{2d}{d-2}$, a contradiction.
\end{proof}

\begin{proof}[Proof of Corollary \ref{pet2}]  From Proposition \ref{precompac} we see that (i) implies (iv), which trivially implies (iii).  The implication of (ii) from (iii) follows from \eqref{uplus}, Lemma \ref{rl}, and duality.  So it suffices to show that (ii) implies (i).  Accordingly, let $u$ be a forward-global solution (of energy at most $E$, say) with $u_+ = 0$, then by \eqref{jgk} we have
$$ \lim_{t \to +\infty} \dist_H(u(t), J(GK)) = 0$$
for some compact set $K \subset H$.  Thus there exists an increasing sequence of times $T_n \to \infty$ such that
$$ \dist_H(u(t), J(GK)) \leq 2^{-n} \hbox{ whenever } t > T_n.$$
Now the partial orbit $\{ u(t): 0 \leq t \leq T_n \}$ is compact for each $n$, and so we can
find a compact subset $K_n$ of $J(GK)$ such that
$$ \dist_H(u(t), J(GK)) = \dist_H(u(t), K_n) + O(2^{-n}) \hbox{ whenever } 0 \leq t \leq T_n.$$
We can easily arrange so that the $K_n$ are increasing in $n$.  We can thus split $u(t) = w(t) + y(t)$ for all $t \in \R$, where $w(t) \in K_n$ whenever $t \leq T_n$, and $\|y(t)\|_H \lesssim 2^{-n}$ whenever $t > T_n$.  Since $K_n$ and $\{ u(t): 0 \leq t \leq T_n\}$ are both compact, we easily see that
$\{ y(t): 0 \leq t \leq T_n \}$ is covered by finitely many balls of radius $O(2^{-n})$ for each $n$.  Since $\{ y(t): t > T_n \}$ is covered by a single ball of radius $O(2^{-n})$, we conclude that the orbit $\{ y(t): 0 \leq t < \infty\}$ is totally bounded and hence contained in a compact set $K' \subset H$.  We then have
$$ \{ u(t): 0 \leq t < \infty \} \subset J(GK) + K' \subset (J+1)( G (K \cup K') )$$
and so $u$ is almost periodic as claimed.
\end{proof}

\begin{proof}[Proof of Corollary \ref{pet1}]  From Proposition \ref{precompac} we see that (i) implies (iii), which trivially implies (ii).  Suppose now that (ii) holds.  From Corollary \ref{pet2} we already have that $u_+ = 0$ and $u$ is $G$-almost periodic and thus lies in $J(GK)$ for some $J \geq 1$ and compact $K$.  Applying Proposition \ref{precompac} to $K$ we conclude that for every $\mu_0 > 0$ there exists $\mu_1 > 0$ such that
$$ \| P_{\geq 1/\mu_1} f \|_H \lesssim \mu_0$$
for all $f \in K$, and hence by the triangle inequality
$$ \| P_{\geq 1/\mu_1} u(t) \|_H \lesssim J\mu_0$$
for all $t \geq 0$.  Applying Proposition \ref{precompac} one last time and using (ii) we obtain that the orbit of $u$ is precompact, as desired.
\end{proof}

\appendix

\section{Dispersive estimates}\label{dispersive}

In this section we recall some standard dispersive estimates for the Schr\"odinger equation.
We begin by recalling the standard fixed-time estimates
\begin{equation}\label{energy-dispersive}
\| e^{it\Delta} f \|_{L^{r'}_x(\R^d)} \lesssim \frac{1}{|t|^{d(\frac{1}{r}-\frac{1}{2})}} \| f \|_{L^r_x(\R^5)} \hbox{ for all }
1 \leq r \leq 2
\end{equation}
where $r' := \frac{r}{r-1}$ is the dual exponent to $r$.  Indeed, the case $r=2$ follows from \eqref{prop-fourier} and Plancherel's theorem, while the case $r=1$ follows from \eqref{prop}, and the intermediate cases then follow by interpolation.

We recall that the fixed time estimates \eqref{energy-dispersive} imply the \emph{Strichartz estimates}
\begin{align}
 \| e^{it\Delta} f \|_{L^q_t L^r_x(\R \times \R^d)} &\lesssim \|f\|_{L^2_x(\R^d)}\label{homog}\\
 \| \int_\R e^{-it\Delta} F(t) \|_{L^2(\R^d)} &\lesssim \|F\|_{L^{\tilde q'}_t L^{\tilde r'}_x(\R \times \R^d)} \label{dual-homog} \\
 \| \int_{t' < t} e^{i(t-t')\Delta} F(t') \|_{L^q_t L^r_x(\R \times \R^d)} &\lesssim 
\|F\|_{L^{\tilde q'}_t L^{\tilde r'}_x(\R \times \R^d)}; \label{retarded}
\end{align}
whenever $(q,r)$ and $(\tilde q, \tilde r)$ are admissible and for any test functions (say) $f, F$; see e.g. \cite{tao:keel} for a proof.  

It is now well-known to the experts that the linear Strichartz estimates come with bilinear refinements,
which roughly speaking assert that the ``high-low'' interactions of different frequency components are weaker than one might first expect.  There are several such bilinear estimates available; it shall be convenient to use the following version, due to Visan \cite{visan-thesis}:

\begin{theorem}[Bilinear Strichartz estimate]\label{bilst}  For any interval $I \subset \R$, any $t_0 \in I$, and any $0 < \delta \leq \frac{1}{2}$, we have
\begin{align*}
\| uv \|_{L^2_{t,x}(I \times \R^d)} &\lesssim_{\delta,q,r,\tilde q,\tilde r} 
(\|u(t_0)\|_{\dot H^{-1/2+\delta}_x(\R^d)} + \| |\nabla|^{-1/2+\delta} (i\partial_t+\Delta) u \|_{L^q_t L^r_x(I \times \R^d)}) \\
&\quad \times
(\|v(t_0)\|_{\dot H^{(d-1)/2-\delta}_x(\R^d)} + \| |\nabla|^{(d-1)/2 - \delta} (i\partial_t+\Delta) v \|_{L^{\tilde q}_t L^{\tilde r}_x(I \times \R^d)}) 
\end{align*}
for any $u,v$ and any admissible pairs $(q,r)$, $(\tilde q, \tilde r)$ with $q, \tilde q > 2$.
Here $\|f\|_{\dot H^s_x} := \| |\nabla|^s f\|_{L^2_x}$ denotes the homogeneous Sobolev norms.
\end{theorem}

\begin{proof} See \cite[Lemma 2.5]{visan-thesis}; the proof combines a standard bilinear Strichartz estimate (see \cite{borg:book}, \cite{gopher}) with the Christ-Kiselev lemma.
\end{proof}

As it turns out we will be content with the (relatively weak) $\delta=1/2$ case of this theorem.

\subsection{Perturbation theory}

We now use the Strichartz estimates to establish some standard perturbation theory results for
the NLS equation \eqref{nls}.

\begin{lemma}[Perturbation lemma]\label{perturb}  Let $u_0 \in H$, let $I$ be a compact time interval containing a time $t_0$, let $A, \mu_0, \mu_1 > 0$, and let $v: I \to H$ be a strong solution to the forced NLS
$$ iv_t + \Delta v = F(v) + G$$
where $v, G$ obey the estimates
\begin{equation}\label{va}
 \| v \|_{L^{q_0}_t W^{1,r_0}_x \cap C^0_t H^1_x( I \times \R^d )} + \|v(t_0)-u_0\|_H 
\lesssim A
\end{equation}
and
\begin{equation}\label{qeps}
 \|e^{i(t-t_0)\Delta} (v(t_0)-u_0)\|_{L^{q_0}_t W^{1,r_0}_x \cap L^\infty_t L^Q_x(I \times \R^d)}
 \lesssim \mu_1
\end{equation}
and 
\begin{equation}\label{geps}
\| G \|_{L^{q'_0}_t W^{1,r'_0}_x(I \times \R^d)} 
\lesssim \mu_1.
\end{equation}
Then, if $\mu_1$ is sufficiently small depending on $A$, $|I|$, $\mu_0$, there exists a solution $u: I \to H$ to \eqref{nls} with $u(t_0)=v(t_0)$ and
$$ \| u-v \|_{L^{q_0}_t W^{1,r_0}_x \cap C^0_t H^1_x( I \times \R^d )} \lesssim \mu_0.$$
\end{lemma}

\begin{remark} In view of \eqref{homog} and Sobolev embedding, the hypothesis \eqref{qeps} is a consequence of the simpler hypothesis $\|v(t_0) - u_0\|_H \lesssim \mu_1$.  In the case of critical NLS, this type of perturbation argument appeared in \cite{TV}; the subcritical case is easier, though there are technical difficulties arising from the fact that the exponent $\theta$ in \eqref{fpower-3} can be strictly less than one (which in particular must occur in the subquadratic case $p<2$).
\end{remark}

\begin{proof}  A standard iteration argument (chopping $I$ up into smaller intervals, see e.g. \cite{gopher}, \cite{visan-thesis}, \cite{TV}) shows that it suffices to establish this estimate assuming that $|I|$ is sufficiently small depending on $A$.  In particular the local theory (Theorem \ref{local}) now ensures that the solution $u$ exists on $I$.

We make the ansatz $u = v+w$, then $w$ solves the equation
$$ iw_t + \Delta w = F(v+w) - F(v) - G; \quad w(t_0) = v(t_0)-u_0.$$
Now introduce the quantity
$$ X := \| w \|_{L^{q_0}_t W^{1,r_0}_x \cap L^\infty_t H^1_x(I \times \R^d)}.$$
From the local theory we have
\begin{equation}\label{xic}
 X \lesssim_A 1.
\end{equation}

We would like to improve this bound on $X$, but first we must control the lower-order quantity
$$ Y := \| w \|_{L^{q_0}_t L^{r_0}_x(I \times \R^d)}.$$
By the Strichartz estimate \eqref{retarded} and \eqref{qeps}, \eqref{geps}, we have
$$ Y \lesssim \mu_1 + \| F(v+w) - F(v) \|_{L^{q'_0}_t L^{r'_0}_x(I \times \R^d)}.$$
From \eqref{fpower}, \eqref{fpower-2} we can bound
$$ |F(v+w)-F(v)| \lesssim (|v|+|w|)^{p-1} |w|$$
and hence by H\"older's inequality
$$ Y \lesssim_A \mu_1 + |I|^{1/q_0 - 1/q'_0} (1 + X)^{p-1} Y.$$
Using \eqref{xic} and a continuity argument, we conclude that
\begin{equation}\label{ya}
Y \lesssim_A \mu_1.
\end{equation}

Now we revisit $X$.  Applying the Strichartz estimate \eqref{retarded} and \eqref{qeps}, \eqref{geps} we conclude that 
$$ X \lesssim \mu_1 + \| F(v+w) - F(v) \|_{L^{q'_0}_t W^{1,r'_0}_x(I \times \R^d)}.$$
On the other hand, from \eqref{fpower}, \eqref{fpower-2}, \eqref{fpower-3} we have the pointwise bound
$$ |\nabla^j (F(v+w) - F(v))| \lesssim (|v|+|w|)^{p-1} |\nabla^j w| + (|v|+|w|)^{p-1-\theta} |w|^\theta |\nabla^j v|$$
for $j=0,1$.  Using H\"older's inequality and \eqref{xic}, \eqref{va} we can bound
$$ \| (|v|+|w|)^{p-1} |\nabla^j w| \|_{L^{q'_0}_t L^{r'_0}_x(I \times \R^d)} \lesssim_A |I|^{1/q_0-1/q'_0} X$$
while from H\"older's inequality and \eqref{xic}, \eqref{va} (and \eqref{admis}) we can bound
$$ \| (|v|+|w|)^{p-1-\theta} |w|^\theta |\nabla^j v| \|_{L^{q'_0}_t L^{r'_0}_x(I \times \R^d)} \lesssim_A \| w \|_{L^s_t L^Q_x(I \times \R^d)}^\theta $$
for some $s<\infty$ depending only on $p, \theta, q_0, r_0$.  However, from interpolation between \eqref{xic}, \eqref{ya} and Sobolev embedding we see that
$$ \| w \|_{L^s_t L^Q_x(I \times \R^d)} \lesssim_A \mu_1^c$$
for some $c>0$ depending only on $p, \theta, q_0, r_0$.  We thus have
$$ X \lesssim_A \mu_1 + |I|^{1/q_0-1/q'_0} X + \mu_1^c$$
which (for $|I|$ sufficiently small depending on $A$) gives the bound
$$ X \lesssim_A \mu_1 + \mu_1^c$$
and hence $X \leq \mu_0$ if $\mu_1$ is sufficiently small depending on $A, \mu_0$.  This gives the claim.
\end{proof}

\section{Compact subsets of $H$}\label{appendix}

The purpose of this appendix is to establish some basic properties of precompact, $G$-precompact, and compact subsets of the energy space $H$.  We first show that precompactness in $H$ is equivalent to simultaneous localisation in both space and frequency (cf. \cite{compact}).

\begin{proposition}[Equivalence of precompactness and localisation]\label{precompac}  Let $K \subset H$.  Then the following are equivalent:
\begin{itemize}
\item[(i)] $K$ is precompact in $H$ (i.e. $K$ is contained in a compact subset of $H$).
\item[(ii)] $K$ is bounded, and for any $\mu_0 > 0$ there exists $\mu_1 > 0$ such that we have the frequency localisation estimate
$$ \| P_{\geq 1/\mu_1} f \|_H \lesssim \mu_0$$
and the spatial localisation estimate
$$ \int_{|x| \geq 1/\mu_1} |f(x)|^2\ dx \lesssim \mu_0^2$$
for all $f \in K$.
\item[(iii)] $K$ is bounded, and for any $\mu_0 > 0$ there exists $\mu_1 > 0$ such that we have the frequency localisation estimates
$$ \| P_{\geq 1/\mu_1} f \|_H \lesssim \mu_0$$
and
$$ \| P_{\leq \mu_1} f \|_H \lesssim \mu_0$$
and the improved spatial localisation estimate
$$ \int_{|x| \geq 1/\mu_1} |f(x)|^2 + |\nabla f(x)|^2\ dx \lesssim \mu_0^2$$
for all $f \in K$.
\end{itemize}
\end{proposition}

\begin{proof} Let us first show that (i) implies (iii).  A simple application of the monotone convergence theorem shows that (iii) holds when $K$ is a singleton set, and hence when $K$ is finite.  From the triangle inequality we conclude that for any fixed $\mu_0 > 0$, that (iii) holds whenever $K$ is covered by finitely many balls of radius $\mu_0$.  Since precompact sets are totally bounded, the claim follows.

Clearly (iii) implies (ii), so it remains to show that (ii) implies (i).  Assume $K$ is such that (ii) holds.  It suffices to show that $K$ is totally bounded, and in particular it will suffice to show that $K$ can be covered by finitely many balls of radius $\mu_0$ for any fixed $\mu_0 > 0$.  By (ii) we can find $\mu_1 > 0$ such that
\begin{equation}\label{com1}
 \| P_{\geq 1/\mu_1} f \|_H \lesssim \mu_0
\end{equation}
for all $f \in K$, and then by (ii) again we can find $\mu_2 > 0$ such that
\begin{equation}\label{com2}
\int_{|x| \geq 1/\mu_2} |f(x)|^2\ dx \lesssim \mu_0^2 \mu_1^2
\end{equation}
for all $f \in K$.  Let $\chi$ be a bump function supported on the ball $\{ |x| \leq 2/\mu_2\}$ which equals one on $|x| \leq 1/\mu_2$.  Given any $f \in K$, we split
$$f = P_{< 1/\mu_1}[\chi f] + P_{<1/\mu_1}[(1-\chi) f] + P_{\geq 1/\mu_1} f.$$
By \eqref{com1}, the third term is $O_H(\mu_0)$; by \eqref{com2}, the second term is also 
$O_H(\frac{1}{\mu_1} \mu_0 \mu_1) = O_H(\mu_0)$.  Finally, from Rellich embedding (or the Arzela-Ascoli theorem) one easily verifies that $P_{<1/\mu_1} \chi$ is a compact operator on $H$, and so (as $K$ is bounded) the set $P_{<1/\mu_1} \chi K$ is covered by finitely many balls of radius $O(\mu_0)$.  The claim follows.
\end{proof}

We shall actually need a generalisation of the above proposition, in which the individual functions $f$ are replaced by trajectories $u(t)$ in $H$.

\begin{proposition}[Criterion for compact attractor]\label{attractive}  Let ${\mathcal U}$ be a collection of trajectories $u: \R^+ \to H$.  Then the following are equivalent:
\begin{itemize}
\item[(i)] There exists a compact set $K \subset H$ such that $\lim_{t \to +\infty} \dist_H(u(t),K) = 0$ for all $u \in {\mathcal U}$.
\item[(ii)] ${\mathcal U}$ is asymptotically bounded in the sense that
\begin{equation}\label{asymb}
\sup_{u \in {\mathcal U}} \limsup_{t \to +\infty} \| u(t) \|_H < \infty
\end{equation}
and for any $\mu_0 > 0$ there exists $\mu_1 > 0$ such that we have the asymptotic frequency localisation estimate
$$ \limsup_{t \to +\infty} \| P_{\geq 1/\mu_1} u(t) \|_H \lesssim \mu_0$$
and the spatial localisation estimate
$$ \limsup_{t \to +\infty} \int_{|x| \geq 1/\mu_1} |u(t,x)|^2\ dx \lesssim \mu_0^2$$
for all $u \in {\mathcal U}$.
\item[(iii)] ${\mathcal U}$ is asymptotically bounded in the sense of \eqref{asymb},
and for any $\mu_0 > 0$ there exists $\mu_1 > 0$ such that we have the asymptotic frequency localisation estimates
$$ \limsup_{t \to +\infty} \| P_{\geq 1/\mu_1} u(t) \|_H \lesssim \mu_0$$
and 
$$ \limsup_{t \to +\infty} \| P_{\leq \mu_1} u(t) \|_H \lesssim \mu_0$$
the improved spatial localisation estimate
$$ \limsup_{t \to +\infty} \int_{|x| \geq 1/\mu_1} |u(t,x)|^2 + |\nabla u(t,x)|^2 
\ dx \lesssim \mu_0^2$$
for all $u \in {\mathcal U}$.
\end{itemize}
\end{proposition}

Note that Proposition \ref{precompac} is essentially the special case of Proposition \ref{attractive} when the $u(t)$ are constant in $t$.

\begin{proof} The implication of (iii) from (i) follows immediately from the analogous implication in
Proposition \ref{precompac}, while the implication of (ii) from (iii) is still trivial.  Thus it only remains to deduce (i) from (ii).

Let ${\mathcal U}$ be such that (ii) holds.  By repeating the arguments in Proposition \ref{precompac}
we see that for any $\mu_0 > 0$ there exists a subset $E_{\mu_0}$ of $H$ which is covered by a finite union of balls of radius $\mu_0$, which absorbs ${\mathcal U}$ in the sense that for all $u \in {\mathcal U}$ we have $u(t) \in E_{\mu_0}$ for all sufficiently large $t$.

Now let $F_n := \bigcap_{1 \leq m \leq n} \overline{E_{2^{-m}}}$ for each integer $n > 0$, then the $F_n$ are a nested sequence of closed subsets of $H$, with each $F_n$ covered by finitely many balls of radius $O(2^{-n})$.  If we let $K := \bigcap_n F_n$, then we conclude that $K$ is closed and totally bounded, hence compact. Furthermore, we see that any
sequence $f_n \in F_n$ is also totally bounded and hence every subsequence has a further convergent subsequence, whose limit must necessarily lie in $K$. Taking contrapositives, we conclude that every open neighbourhood of $K$ must contain $F_n$ (and hence $E_{2^{-n}}$) for sufficiently large $n$.  From this we easily see that $\lim_{t \to +\infty} \dist(u(t),K) = 0$ for all $u \in {\mathcal U}$ as desired.
\end{proof}

Finally, we need a $G$-precompact analogue of the above proposition:

\begin{proposition}[Criterion for $G$-compact attractor]\label{G-attractive}  Let ${\mathcal U}$ be a collection of trajectories $u: \R^+ \to H$, and let $J \geq 1$.  Then the following are equivalent:
\begin{itemize}
\item[(i)] There exists a $G$-precompact set $K \subset H$ with $J$ components such that $\lim_{t \to +\infty} \dist_H(u(t),K) = 0$ for all $u \in {\mathcal U}$.
\item[(ii)] ${\mathcal U}$ is asymptotically bounded in the sense of \eqref{asymb},
such that for any $\mu_0 > 0$ there exists $\mu_1 > 0$ such that 
for every $u \in {\mathcal U}$ we have $x_1,\ldots,x_J: \R^+ \to \R^d$ for which we have the asymptotic frequency localisation estimate
\begin{equation}\label{lsu}
 \limsup_{t \to +\infty} \| P_{\geq 1/\mu_1} u(t) \|_H \lesssim \mu_0
 \end{equation}
and the spatial localisation estimate
\begin{equation}\label{ssu}
\limsup_{t \to +\infty} \int_{\inf_{1 \leq j \leq J} |x-x_j(t)| \geq 1/\mu_1} |u(t,x)|^2\ dx \lesssim \mu_0^2.
\end{equation}
\item[(iii)] ${\mathcal U}$ is asymptotically bounded in the sense of \eqref{asymb},
such that for any $\mu_0 > 0$ there exists $\mu_1 > 0$ such that 
for every $u \in {\mathcal U}$ we have $x_1,\ldots,x_J: \R^+ \to \R^d$ for which we have the asymptotic frequency localisation estimates
$$ \limsup_{t \to +\infty} \| P_{\geq 1/\mu_1} u(t) \|_H \lesssim \mu_0$$
and
$$ \limsup_{t \to +\infty} \| P_{\leq \mu_1} u(t) \|_H \lesssim \mu_0$$
and the improved spatial localisation estimate
$$ \limsup_{t \to +\infty} \int_{\inf_{1 \leq j \leq J} |x-x_j(t)| \geq 1/\mu_1} 
|u(t,x)|^2 + |\nabla u(t,x)|^2\ dx \lesssim \mu_0^2.$$
\end{itemize}
\end{proposition}

\begin{proof}  The implication of (iii) from (i) follows from the analogous implications from preceding propositions, together with the triangle inequality.  As the implication of (ii) from (iii) is still trivial, it once again suffices to show that (ii) implies (i).

Let $u \in {\mathcal U}$, and let $x_1,\ldots,x_J: \R^+ \to \R^d$ be as in (ii).  We form the partition of unity
$$ 1 = \sum_{j=1}^J \psi_{j,t}(x)$$
where
$$ \psi_{j,t}(x) := \frac{ \langle x - x_j(t) \rangle^{-1} }{\sum_{j'=1}^J \langle x - x_{j'}(t) \rangle^{-1} }.$$
Clearly $\psi_{j,t}$ ranges between zero and one.  Direct computation also yields the $C^2$ regularity bounds
\begin{equation}\label{nablak}
 \nabla^k \psi_{j,t} = O_J(1) \hbox{ for } k=0,1,2.
\end{equation}

We can now split
\begin{equation}\label{vdecomp}
 u(t) = \sum_{j=1}^J \tau_{x_j(t)} w_j(t)
 \end{equation}
where
$$ w_j(t) := \tau_{-x_j(t)}( \psi_{j,t} u(t) ).$$

\begin{lemma}[Localisation of $w_j(t)$]\label{wjt}  For any $\mu_0 > 0$ there exists $\mu_2 > 0$ such that
\begin{equation}\label{pwj}
 \limsup_{t \to +\infty} \| P_{\geq 1/\mu_2} w_j(t) \|_H \lesssim \mu_0
 \end{equation}
and
\begin{equation}\label{wj-spat}
 \limsup_{t \to +\infty} \int_{|x| \geq 1/\mu_2} |w_j(t,x)|^2\ dx \lesssim \mu_0^2
 \end{equation}
for all $1 \leq j \leq J$ and all $u \in {\mathcal U}$, where $w_j$ is defined as above.
\end{lemma}

\begin{proof}  Fix $\mu_0$, and choose $\mu_1$ sufficiently small depending on $\mu_0$, and then $\mu_2$ sufficiently small depending on $\mu_0,\mu_1$.  Let $u \in {\mathcal U}$, and let $J, x_j, \psi_{j,t}, w_j$ be as above, and fix $1 \leq j \leq J$.  By \eqref{lsu} (choosing $\mu_2$ small enough) we have
$$ \| P_{\geq 1/100\mu_2} u(t) \|_H \lesssim \mu_0$$
for all sufficiently late times $t$. By \eqref{nablak} and the Leibnitz rule this ensures that
$$ \| P_{\geq 1/\mu_2} ( \psi_{j,t} P_{\geq 1/100\mu_2} u(t)) \|_H \lesssim \mu_0.$$
Now consider the expression
$$ \| P_{\geq 1/\mu_2} ( \psi_{j,t} P_{< 1/100\mu_2} u(t)) \|_H.$$
By Fourier analysis we see that we may freely replace $\psi_{j,t}$ by $P_{\geq 1/10\mu_2} \psi_{j,t}$.  If we then discard the bounded multiplier $P_{\geq 1/\mu_2}$ and using the Leibnitz rule we can bound
$$ \| P_{\geq 1/\mu_2} ( \psi_{j,t} P_{< 1/100\mu_2} u(t)) \|_H
\lesssim \| P_{\geq 1/10\mu_2} \psi_{j,t} \|_{C^1_x(\R^d)} \|u(t)\|_H.$$
By \eqref{nablak} and \eqref{asymb} we conclude that
$$ \| P_{\geq 1/\mu_2} ( \psi_{j,t} P_{< 1/100\mu_2} u(t)) \|_H \lesssim \mu_2$$
for all sufficiently late times $t$. By the triangle inequality we conclude
$$ \| P_{\geq 1/\mu_2} ( \psi_{j,t} u(t)) \|_H \lesssim \mu_0;$$
translating this by $-x_j(t)$ we obtain \eqref{pwj}.  

Now we prove \eqref{wj-spat}.  Translating by $x_j(t)$, it suffices to show that
$$ \int_{|x - x_j(t)| \geq 1/\mu_2} \psi_{j,t}(x)^2 |u(t,x)|^2\ dx \lesssim \mu_0^2$$
for sufficiently late times $t$.  Let $D := \inf_{1 \leq j' \leq J} |x-x_{j'}(t)|$.  From
\eqref{ssu} we have (if $\mu_1$ is sufficiently small depending on $\mu_0$)
$$
\int_{\R^d} 1_{D > 1/\mu_1} |u(t,x)|^2\ dx \lesssim \mu_0^2
$$
for sufficiently late times $t$.
Since $\psi_{j,t}(x) = O(1)$, we thus see from the triangle inequality that it suffices to show that
$$ \int_{|x - x_j(t)| \geq 1/\mu_2} \psi_{j,t}(x)^2 1_{D \leq 1/\mu_1} |u(t,x)|^2\ dx \lesssim \mu_0^2.$$
However, by construction of $\psi_{j,t}$ we see that 
$\psi_{j,t}(x) = O( \mu_2 / \mu_1 )$ whenever $D \leq 1/\mu_1$ and $|x-x_j(t)| \geq 1/\mu_2$.  The claim now follows from \eqref{asymb}.
\end{proof}

From the above lemma and Proposition \ref{attractive},
we conclude the existence of a compact set $K \subset H$ such that
$$ \lim_{t \to +\infty} \dist_H(w_j(t), K) = 0$$
for all $1 \leq j \leq J$ and all forward-global solutions $u$ of energy at most $E$.  From this and
\eqref{vdecomp} we see that
$$ \lim_{t \to +\infty} \dist_H(v(t), J(GK)) = 0.$$
Since $J(GK)$ is $G$-precompact with $J$ components, the claim follows.
\end{proof}

Next, we show that the linear propagator $e^{it\Delta}$, when applied to precompact sets, sends them to zero in weak norms:

\begin{lemma}[Riemann-Lebesgue lemma]\label{rl}  Let $K$ be a precompact subset of $H$.  Then for any $2 < q \leq \frac{2d}{d-2}$ and $R > 0$, we have
$$ \lim_{t \to \pm \infty} \sup_{f \in K} \| e^{it\Delta} f \|_{L^q(\R^d)} = 0$$
and
$$ \lim_{t \to \pm \infty} \sup_{f \in K} \sup_{x_0 \in \R^d} \int_{|x-x_0| \leq R} |e^{it\Delta} f(x)|^2 + |\nabla e^{it\Delta} f(x)|^2\ dx = 0.$$
\end{lemma}

\begin{proof} Since $e^{it\Delta}$ is unitary, we see from Sobolev embedding that
$$ \| e^{it\Delta} f \|_{L^q(\R^d)} \lesssim \|f\|_H$$
and
$$ \int_{|x-x_0| \leq R} |e^{it\Delta} f(x)|^2 + |\nabla e^{it\Delta} f(x)|^2\ dx \lesssim \|f\|_H^2.$$
Since precompact sets are totally bounded, a standard argument then shows that to establish the claims it suffices to do so for finite sets $K$, and hence for singleton sets $K$.  By a limiting argument it
then suffices to verify the claim when $K = \{f\}$ and $f \in C^\infty_0(\R^d)$ is a test function.  But then the claim follows from direct computation using \eqref{prop} and stationary phase (one can also
use \eqref{energy-dispersive}).
\end{proof}

By the translation invariance of the quantities in the above lemma, and the triangle inequality, we also conclude

\begin{corollary}[Riemann-Lebesgue lemma for $G$-precompact sets]\label{grl}  Let $K$ be a $G$-precompact subset of $H$ with $J$ components.  Then for any $2 < q \leq \frac{2d}{d-2}$ and $R > 0$, we have
$$ \lim_{t \to \pm \infty} \sup_{f \in K} \| e^{it\Delta} f \|_{L^q(\R^d)} = 0$$
and
$$ \lim_{t \to \pm \infty} \sup_{f \in K} \sup_{x_0 \in \R^d} \int_{|x-x_0| \leq R} |e^{it\Delta} f(x)|^2 + |\nabla e^{it\Delta} f(x)|^2\ dx = 0.$$
\end{corollary}

Next, we analyse the convergence properties of sequences in $G$-precompact sets.

\begin{lemma}[Baby concentration compactness]\label{cc}  Let $K$ be a compact subset of $H$, let $J \geq 1$, and let $f_n$ be a sequence in $J(GK)$.  Then, after passing to a subsequence, there exists a partition $J = J_1 + \ldots + J_M$ with $J_1,\ldots,J_M \geq 1$, functions $w_m \in J_m(GK)$ and points
$x_{m,n} \in \R^d$ for $1 \leq m \leq M$ and $n \geq 1$ such that we have the decomposition
\begin{equation}\label{fn}
 f_n = \sum_{m=1}^M \tau_{x_{m,n}} w_m + o_H(1)
 \end{equation}
where the error $o_H(1)$ goes to zero in $H$ norm as $n \to \infty$,
and such that we have the asymptotic separation condition
\begin{equation}\label{ass}
\lim_{n \to \infty} |x_{m,n} - x_{m',n}| \to \infty
\end{equation}
for all $1 \leq m < m' \leq M$.
\end{lemma}

\begin{proof}  By definition of $J(GK)$, we have a representation
$$ f_n = \sum_{j=1}^J \tau_{y_{j,n}} u_{j,n}$$
for some $y_{j,n} \in \R^d$ and $u_{j,n} \in K$ for $j=1,\ldots,J$.  By passing to a subsequence repeatedly and exploiting the compactness of $K$ we may assume for each $1 \leq j \leq J$ that $u_{j,n} \to w_j$ as $n \to \infty$ for some $u_j \in K$.  Since the contribution of the error $u_{j,n}-u_j$ is $o(1)$, we may thus assume that
$$ f_n = \sum_{j=1}^J \tau_{y_{j,n}} u_j.$$
Next, by exploiting the local compactness of $\R^d$ and repeatedly passing to subsequences we may assume for each $1 \leq j,j' \leq J$ that either $y_{j,n} - y_{j',n}$ converges to an element of $\R^d$, or else diverges to infinity.  The former case determines an equivalence relation on $\{1,\ldots,J\}$, and let $\{1,\ldots,J\} = A_1 \cup \ldots \cup A_M$ be the associated equivalence classes.  We may then write
$$ y_{j,n} = x_{m,n} + z_j + o(1)$$
for all $n \geq 1$, $1 \leq m \leq M$, and $j \in A_m$, where $x_{m,n}, z_j \in \R$ and $o(1)$ goes to zero as $n \to \infty$.  Again, we may absorb the effect of the $o(1)$ error into the $o_H(1)$ term of
\eqref{fn} (because the action of the translation group is continuous in the strong operator topology).
We thus have
$$ f_n = \sum_{m=1}^M \tau_{x_{m,n}} u_m$$
where $w_m := \sum_{j \in A_m} \tau_{z_j} u_j$.  Setting $J_m := |A_m|$ we see that $w_m \in J_m(GK)$
and $J_1 + \ldots + J_M = J$.  The asymptotic separation condition \eqref{ass} then follows from the definition of the equivalence relation, and the claim follows.
\end{proof}

\begin{corollary}[Compact sets have closed multi-orbits]\label{compact-k}  Let $K$ be a compact subset of $H$.  Then $J(GK)$ is closed for all $J \geq 1$.
\end{corollary}

\begin{proof}  Let $f_n$ be a sequence in $J(GK)$ which converges to some $f \in H$, thus $f_n = f + o_H(1)$.  By passing to a subsequence, we may invoke Lemma \ref{cc} and obtain a decomposition \eqref{fn}.  By passing to a further subsequence, we may assume that each sequence $(x_{m,n})_{n \geq 1}$ is either convergent to some limit $x_m$, or else goes to infinity.  The first condition can only occur for at most one $m$.  By permuting the $m$ and adding a dummy index if necessary we may assume that it occurs for $m=1$.  By discarding a $o(1)$ error as before we may assume that the $x_{1,n}$ are in fact constant in $n$, and then by absorbing this constant into $w_1$ we may take $x_{1,n} = 0$.  Thus
we have
$$ w_1 - f + \sum_{m=2}^M \tau_{x_{m,n}} w_m = o_H(1).$$
Taking $L^2$ norms and using the asymptotic separation \eqref{ass} we see that
$$ \lim_{n \to \infty} \| w_1 - f + \sum_{m=2}^M \tau_{x_{m,n}} w_m \|_H^2 - 
\| w_1 - f \|_H^2 - \sum_{m=2}^M \| w_m \|_H^2 = 0$$
and we therefore conclude that $w_1 = f$ and $w_m=0$ for $m \geq 2$.  In particular $f = w_1 + \ldots + w_m \in J(GK)$ and the claim follows.
\end{proof}

\begin{corollary}\label{gkclosed} The closure of any $G$-precompact set with $J$ components is also $G$-precompact with $J$ components.
\end{corollary}

We also need a kind of converse to Lemma \ref{cc}:

\begin{lemma}[Asymptotic profiles of $J(GK)$]\label{cc-inverse} Let $K$ be a compact subset of $H$, let $J \geq 1$, and let $f_n$ be a sequence in $J(GK)$ which has a decomposition \eqref{fn}, where $x_{m,n}$ obeys \eqref{ass} and $w_m \in H$.  Then there exists a partition $J = J_1 + \ldots + J_M$ with $J_1,\ldots,J_M \geq 0$ such that $w_m \in J_m(GK)$.
\end{lemma}

\begin{proof}  We apply Lemma \ref{cc} and (after passing to a subsequence) we obtain an alternate decomposition
\begin{equation}\label{fnp}
 f_n = \sum_{m'=1}^{M'} \tau_{y_{m',n}} u_{m'} + o_H(1)
\end{equation}
for some partition $J = J'_1 + \ldots + J'_{M'}$, where $u_{m'} \in J'_{m'}(GK)$ and $|y_{m',n} - y_{m''_n}| \to \infty$ as $n \to \infty$ for all $1 \leq m' < m'' \leq M'$.

By passing to further subsequences we may assume that for each $1 \leq m \leq M$ and $1 \leq m' \leq M'$, we either have
$x_{m,n}-y_{m',n}$ converge in $\R^d$, or else diverge to infinity.  Suppose that we can find $1 \leq m \leq M$ such that $x_{m,n}-y_{m',n}$ diverged for all $m'$.  Then from \eqref{fnp} we see that
$$ \langle f_n, \tau_{x_{m,n}} w_m \rangle_H \to 0 \hbox{ as } n \to \infty.$$
On the other hand, from \eqref{fn} we see that
$$ \langle f_n, \tau_{x_{m,n}} w_m \rangle_H \to \|w_m\|_H^2 \hbox{ as } n \to \infty.$$
Thus $w_m=0$ in this case.  We can discard these cases (by setting $J_m=0$) and so we see that for each $1 \leq m \leq M$ there is at least one $1 \leq m' \leq M'$ such that $x_{m,n}-y_{m',n}$ converged; in fact, since the $y_{m',n}$ are asymptotically separated, there is exactly one such $m'$, and $M' \geq M$.  By relabeling we may take $m'=m$, thus $y_{m,n}-x_{m,n}$ converges for all $1 \leq m \leq M$.  By perturbing the $y_{m,n}$ by $o(1)$ and absorbing the resulting error in \eqref{fnp} into the $o_H(1)$ term we may assume that $x_{m,n}-y_{m,n}$ is constant in $n$; by absorbing this constant into $u_m$ (noting that $J'_m(GK)$ is translation-invariant) we can take this constant to be zero, thus $x_{m,n}=y_{m,n}$ for all $1 \leq m \leq M$.

If $m' > M$, then arguing as before we see from \eqref{fn} that
$$ \langle f_n, \tau_{y_{m',n}} u_{m'} \rangle_H \to 0 \hbox{ as } n \to \infty$$
while from \eqref{fnp} we have
$$ \langle f_n, \tau_{y_{m',n}} u_{m'} \rangle_H \to \|u_{m'}\|_H^2 \hbox{ as } n \to \infty$$
and so $u_{m'} = 0$ for all $m' > M$.  Thus \eqref{fnp} becomes
$$ f_n = \sum_{m=1}^M \tau_{x_{m,n}} u_m + o_H(1).$$
Subtracting this from \eqref{fn} we conclude
$$ \sum_{m=1}^M \tau_{x_{m,n}} (w_m-u_m) = o_H(1).$$
Taking $H$ norms of both sides and using asymptotic separation again we conclude that
$$ \sum_{m=1}^M \| w_m-u_m\|_H^2 = 0$$
and so $w_m = u_m$ for all $1 \leq m \leq M$.  In particular $w_m \in J'_m(GK)$ for $1 \leq m < M$.  Since
$w_M =u_M = u_M + u_{M+1} + \ldots + u_{M'}$, we also have 
$$ w_M \in (J'_M + \ldots + J'_{M'})(GK),$$
and the claim follows.
\end{proof}

Now we show that the nonlinear flow maps $S(t)$ are asymptotically additive with respect to asymptotically separated superpositions.

\begin{lemma}[Asymptotic additivity]\label{asl}  Let $w_1,\ldots,w_M \in H$ obey the bound $\|w_1\|_H,\ldots,\|w_M\|_H \leq E$, and for $1 \leq m \leq M$ let $(x_{n,m})_{n \geq 1}$ be a sequence of points in $\R^d$ obeying the separation condition \eqref{ass}.  Let $f_n$ be a sequence in $H$ such that
$$ f_n = \sum_{m=1}^M \tau_{x_{n,m}} w_m + o_H(1).$$
Then if $t$ is a time with $|t|$ sufficiently small depending on $M$ and $E$, we have
$$ S(t) f_n = \sum_{m=1}^M \tau_{x_{n,m}} S(t) w_m + o_H(1).$$
\end{lemma}

\begin{proof}  Let $I$ be a sufficiently small open time interval containing $0$ (depending on $M$ and $E$).  For $t \in I$, let $u_n(t) := S(t) f_n$ and $v_n(t) := \sum_{m=1}^M \tau_{x_{n,m}} S(t) w_m$.  Then we have
$$ \| v_n(0) - u_n(0)\|_H = o(1)$$
and it will suffice to show that
$$ \| v_n - u_n \|_{C^0_t H^1_x(I \times \R^d)} = o(1).$$
From \eqref{homog} we have
$$
 \|e^{i(t-t_0)\Delta} (v_n(0)-u_n(0))\|_{L^{q_0}_t W^{1,r_0}_x \cap L^\infty_t L^Q_x(I \times \R^d)} = o(1)$$
while from the local theory we have
$$
 \| v \|_{L^{q_0}_t W^{1,r_0}_x \cap C^0_t H^1_x( I \times \R^d )} + \|v(t_0)-u_0\|_H \lesssim_{M,E} 1.$$
In view of Lemma \ref{perturb}, it thus suffices to show that the quantity
$$ G_n := (i \partial_t + \Delta) v_n - F(v_n)$$
obeys the bound
$$\| G_n \|_{L^{q'_0}_t L^{r'_0}_x(I \times \R^d)} = o(1).$$
From \eqref{fpower}, \eqref{fpower-2} and induction we have the pointwise estimate
$$ F( z_1 + \ldots + z_M ) = F(z_1) + \ldots + F(z_M) + \sum_{1 \leq m,m' \leq M: m \neq m'} O_M( |z_m| |z_{m'}|^{p-1} )$$
and thus
$$ |G_n| \lesssim_M \sum_{1 \leq m,m' \leq M: m \neq m'} |\tau_{x_{n,m}} S(t) w_m| |\tau_{x_{n,m'}} S(t) w_{m'}|^{p-1}.$$
Thus by the triangle inequality and translation invariance it suffices to show that
$$ \| |\tau_{x_{n,m}-x_{n,m'}} S(t) w_m| |S(t) w_{m'}|^{p-1} \|_{L^{q'_0}_t L^{r'_0}_x(I \times \R^d)} = o(1).$$
If $S(t) w_m$ and $S(t) w_{m'}$ were compactly supported in space then this would follow immediately from the asymptotic separation of $x_{n,m}$ and $x_{n,m'}$.  The general case then follows from a standard limiting argument, noting from local theory that $S(t) w_m$ lies in $L^{q_0}_t L^{r_0}_x$
and $S(t) w_{m'}$ lies in $L^\infty_t L^{Q_0}_x$ on $I \times \R^d$, which by \eqref{admis} and H\"older is enough to justify the limiting argument.
\end{proof}

\end{document}